\documentclass[11pt]{article}
\usepackage{graphicx}
\usepackage{amsfonts}
\usepackage{bbm}
\usepackage{mathrsfs}
\usepackage{mathrsfs,amsmath,amssymb,amsthm}
\usepackage{amssymb}
\usepackage{amsbsy}
\usepackage{color}
\allowdisplaybreaks
\theoremstyle{plain}
\theoremstyle{remark}  \newtheorem{remark}{\noindent\mbox{Remark}}
\theoremstyle{plain}
\theoremstyle{plain}\newtheorem{lemma}{\noindent\mbox{Lemma}}
\theoremstyle{plain} \newtheorem{theorem}{\noindent\mbox{Theorem}}
\theoremstyle{plain}\newtheorem{proposition}{\noindent\mbox{Proposition}}
\theoremstyle{plain}\newtheorem{corollary}{\noindent\mbox{Corollary}}
\theoremstyle{definition} \newtheorem{definition}{\noindent\mbox{Definition}}

 \def\proof{\noindent{\it Proof.~~}}
 \def\qed{\hfill$\Box$\medskip}
 \def\rto{\rightarrow\infty}
 \def\z{\left}
 \def\y{\right}
 \def\no{\nonumber}
  \def\mb{\mathbf}

\begin{document}
 \title{\textbf{ Cutpoints of (1,2) and (2,1) random walks on the lattice of positive half line}}

\author{   Lanlan \uppercase{Tang}$^*$ $\ \ \&$ \ Hua-Ming \uppercase{Wang}$^{*,\dag}$}
\date{}
\maketitle%
 \footnotetext[1]{School of Mathematics and Statistics, Anhui Normal University, Wuhu 241003, China  }
\footnotetext[2]{Email: hmking@ahnu.edu.cn}

\vspace{-.5cm}

\begin{center}
\begin{minipage}[c]{12cm}
\begin{center}\textbf{Abstract}\quad \end{center}

In this paper, we study (1,2) and (2,1) random walks in varying environments on the lattice of positive half line. We assume that the transition probabilities at site $n$ are asymptotically constants as $n\rightarrow\infty.$ For (1,2) random walk, we get some elaborate  asymptotic behaviours of various  escape probabilities and  hitting probabilities of the walk. Such observations and some delicate analysis of continued fractions and the product of nonnegative matrices enable us to give criteria for finiteness of the number of cutpoints of both (1,2) and (2,1) random walks, which generalize   E. Cs\'aki, A. F\"oldes and P. R\'ev\'esz  [J. Theor. Probab. \textbf{23}:
624-638 (2010)] and H.-M. Wang [Markov Processes Relat. Fields \textbf{25}: 125-148 (2019)]. For  near-recurrent random walks, whenever there are infinitely many cutpoints, we also study the asymptotics of the number of cutpoints in $[0,n].$

\vspace{0.2cm}

\vspace{0.2cm}

\end{minipage}
\end{center}

\section{Introduction}\label{s1}

\subsection{Motivation and Backgrounds}
 In this paper, we consider (2,1) and (1,2) random walks in varying environments on the lattice of positive half line. The transition probabilities of the random walks at site $n$  are assumed to be asymptotically constants as $n\rto.$ The notion {\it(1,2) random walk} means that the left-oriented jumps are always of size 1 and the right-oriented jumps are either of  size 1 or of size 2. The notion {\it (2,1) random walk} can be understood similarly.
 Our purpose is to give criteria for finiteness of the number of cutpoints for both (2,1) and (1,2) random walks and study further the asymptotics of the number of cutpoints in $[0,n]$ as $n\rto.$ Roughly speaking, if the walk never returns to $[0,x]$ after its first entry into $[x+1,\infty),$ then $x$ is a cutpoint. Clearly, in the recurrent case, there is no cutpoint on the path of the walk and in the transient case, by intuition if the walk runs to infinity more quickly, there are more cutpoints. For simple high dimensional random walks, the studies of the cutpoints can be founded in \cite{et,jp,l}. For the nearest-neighbor random walk on  $\mathbb Z_+,$ an example of a transient random walk which  has only finitely many cutpoints was given in \cite{jlp}. Later, a sharp criterion was given in \cite{cfrb} for the finiteness of the number of cutpoints. Recently, some extension was made in \cite{w19} for a random walk whose left-oriented jumps are always of size 1 and the right-oriented jumps are always of size 2.  Most recently,  in \cite{lmw}, cutpoints for some  general adapted stochastic processes with bounded increments were studied. Conditions were given to ensure such processes have finitely or infinitely many cutpoints and some partial results for the expectation number of cutpoints in $[0,x]$ were obtained.

 Although cutpoints for more general random walks were studied in \cite{lmw}, we limit ourselves to (1,2) and (2,1) random walks. On one hand, for such simpler models, we can get sufficient and necessary condition for the finiteness of the number of cutpoints. On the other hand,  in \cite{w19}, for the divergent case it is only shown that with a positive probability, (1,2) random walk has infinitely many cutpoints.  To get an almost-sure result for (1,2) random walk, one needs to study deeply the asymptotic behaviours of  escape probabilities and  hitting probabilities of the walk, which are also interesting. In addition, for some near-recurrent random walks, whenever there are infinitely many cutpoints, we also studied the asymptotics of the number of cutpoints lying in $[0,n].$

 Things of (1,2) and (2,1) random walks are more complicated than those of the nearest-neighbor setting, because the escape probability of the walk from an interval is written in terms of the products of nonhomogeneous 2-by-2 nonnegative matrices, whose entries are hard to evaluate. To overcome such difficulties, our idea is to estimate the entries of the products of matrices in terms of the products of tails of some related continued fractions. Then based on some delicate analysis of the continued fraction  and the products of matrices, we can find the asymptotics of the escape probabilities and the hitting probabilities of the random walk, which is crucial to give criteria for the finiteness of the number of cutpoints.

Before presenting the models, we introduce some conventions and notations which will be used in what follows. When we write $A(n)\sim B(n)$, we mean that $A(n)/B(n)\rightarrow 1$ as $n\rto.$ The notation $A(n)=O(B(n))$ means there is a constant $C>0$ such that $|A(n)|<CB(n)$ for all $n$ large enough. For a set $\{\ \},$ we use $\#\{\ \}$ to denote the number of its elements. As usual,  $\mb e_1=(1,0),$ $\mb e_2=(0,1)$ are canonical unit vectors and for a vector $\mb v,$ $\mb v^t$ denotes its transform. We adopt the convention that empty product equals identity and empty sum equals $0$. Finally, throughout the paper, $0<c<\infty$ is a constant whose value may change from line to line.

\subsection{Models and main results}
To introduce precisely the models, we suppose that $q_k,p_{k1},p_{k2},k\ge 2$ are numbers such that $\forall k\ge2, q_k>0,p_{k1}\ge 0,p_{k2}>0$ and  $q_k+p_{k1}+p_{k2}=1.$
Let $X=\{X_k\}_{k\ge 0}$ be a Markov chain on $\mathbb Z_+:=\{0,1,2,...\}$ starting from some $x_0\in\mathbb Z_+,$ with transition probabilities
\begin{align}
&P(X_{k+1}=1|X_k=0)=P(X_{k+1}=2|X_k=1) =1\nonumber,\\
&P(X_{k+1}=n+1|X_k=n)= q_n,\nonumber\\
&P(X_{k+1}=n-1|X_k=n)= p_{n1}, \nonumber\\
&P(X_{k+1}=n-2|X_k=n)= p_{n2}, n\ge 2,k\ge0.\nonumber
\end{align}
Introduce also another Markov chain $Y=\{Y_k\}_{k\ge0}$ on $\mathbb Z_+,$ starting from some $y_0\in \mathbb Z_+,$ with transition probabilities\begin{align*}
 &P(Y_{k+1}=0|Y_k=1)=P(Y_{k+1}=2|Y_k=0)=1,\\
  &P(Y_{k+1}=n-1|Y_k=n)=q_n,\\
  &P(Y_{k+1}=n+1|Y_k=n)= p_{n1}, \\
  &P(Y_{k+1}=n+2|Y_k=n)= p_{n2},n\ge2,k\ge 0.
\end{align*}
Unless otherwise specified we always assume that both $X$ and $Y$ start from $x_0=y_0=2.$ We call the chain $X$ a (2,1) random walk and $Y$ a (1,2) random walk.
For $k\ge 2,$ introduce matrix
\begin{align}
A_k:=\z(\begin{array}{cc}
          a_k & b_k \\
          1 & 0
        \end{array}
\y)\text{ with }a_k:=\frac{p_{k1}+p_{k2}}{q_k},b_k=\frac{p_{k2}}{q_k}.\label{tn}
\end{align}
Let $\varrho_k$ be the spectral radius (the largest eigenvalue) of the matrix $A_k$. Then we can easily compute that
\begin{align}\label{pk}
\varrho_k=\frac{a_k+\sqrt{a_k^2+4b_k}}{2},
\end{align}
which we will work with. We now introduce the following  condition.
\begin{description}
\addtolength{\itemsep}{-0.5em}
  \item[{(C)}] {\it Suppose that $a_k\rightarrow a$ and $b_k\rightarrow b$ for some numbers $a>0$ and $b>0$ as $k\rto.$}
\end{description}

Clearly, under Condition (C),  we have
\begin{align}
A_k\rightarrow A:=\z(\begin{array}{cc}
          a & b \\
          1 & 0
        \end{array}
\y)\no
\end{align} as $k\rto.$
Let \begin{align}
\varrho=\frac{1}{2}{\z(a+\sqrt{a^2+4b}\y)} \text{ and } \sigma=\frac{1}{2}{\z(a-\sqrt{a^2+4b}\y)}\no
\end{align} be the eigenvalues of the matrix  $A.$ Some easy computation shows that
\begin{align}
 -1< \sigma<0.\no
\end{align}

\begin{remark}
  We remark that under Condition (C),
  \begin{align}
       &a+b=1\text{ implies  }\varrho_k\rightarrow\varrho=1 \text{ and }p_{k1}+2p_{k2}-q_k\rightarrow 0,\label{ro}
  \end{align}
   as $k\rto.$ That is to say, the drift of the random walk at site $k$ is asymptotically zero as $k\rto.$ From this point of view, in this case, we usually call such random walk  a {\it near-recurrent} or {\it near-critical} random walk.
\end{remark}

We consider first (1,2) random walk $Y.$ Since $p_{k2}>0,$ standing at site $k,$ with positive probability the walk jumps to $k+2$ in the next step. So it is necessary to consider precisely the escape probabilities and hitting probabilities of the walk.

For $2\le m\le k\le n+1,$ and $j\in\{n,n+1\},$ set
\begin{align}\label{dqa}
  &Q_k^j(m,n)=P(Y\text{ hits }[n,\infty)\text{ before it hits }[0,m] \text{ at }j\big|Y_0=k),
\end{align}
and write
\begin{equation}\label{dqb}
  Q_k(m,n,+):=Q_k^n(m,n)+Q_k^{n+1}(m,n),
\end{equation} which is the probability that starting from $k,$ the chain $Y$ hits $[n,\infty)$ before $[0,m].$

To introduce the hitting probabilities of the chain $Y,$
for $k\ge 0,$ set $L_{k}=\{2k,2k+1\}.$ For $k\geq1$, let
\begin{align}
  T_{k}=\inf\{n\geq 0:Y_{n}\in L_{k}\}\label{dtk}
\end{align}
be the time that $Y$ enters into $L_{k}$ for the first time. Denote
\begin{align}
&h_{k}(1)=P(Y_{T_{k}}=2k),\ h_{k}(2)=P(Y_{T_{k}}=2k+1),k\ge1;\label{dh1}\\
&\eta_{k,m}(1)=P(Y \text{ enters } [m+1,\infty) \text{ at } m+1 |Y_{0}=k),\label{de1}\\
&\eta_{k,m}(2)=P(Y \text{ enters } [m+1,\infty) \text{ at } m+2 |Y_{0}=k),m\ge k\ge 1.\label{de2}
\end{align}
\begin{proposition}\label{hp0}
Under Condition (C), we have \begin{align}
  \lim_{n\rto}\eta_{n,n}(2)=-\sigma,\lim_{n\rto}h_n(2)=-\frac{\sigma}{1-\sigma}\label{lhe}
\end{align}
and for $k\ge1,$ with $\hat a :=a+\frac{(1-\varrho)(\varrho-\sigma)}{\varrho},$ \begin{align}
  \frac{Q_{k+1}^n(k,k+n)}{Q_{k+1}^{n+1}(k,k+n)}\rightarrow-\sigma^{-1},\ Q_{k+n-1}^{k+n+1}(k,k+n)\rightarrow \tau:=
\left\{\begin{array}{ll}-\sigma, &\text{if }\varrho\ge1,\\
\frac{\sqrt{\hat a^2+4b}-\hat a}{2},& \text{if }\varrho<1,\end{array}\right.\label{ql}
\end{align}
uniformly in $k\ge1$ as $n\rto.$
\end{proposition}

\begin{remark}
  Proposition \ref{hp0} is crucial for us to improve the results of \cite{w19}.
\end{remark}

The following proposition gives criteria  for transience of the chain $X$ and $Y.$
\begin{proposition}\label{rt} Suppose that Condition (C) holds. Then we have
{\bf(i)} the chain $X$ is transient if and only if $\sum_{s=2}^{\infty}\varrho_{2}\cdots \varrho_{s}<\infty;$
{\bf(ii)} the chain $Y$ is transient if and only if $\sum_{s=2}^{\infty}\varrho_{2}^{-1}\cdots \varrho_{s}^{-1}<\infty.$
\end{proposition}

Next we study the cutpoints of the chain $X$ and $Y.$ The cutpoints for (2,1) random walk defined below correspond to the strong cutpoints in \cite{cfrb} for nearest-neighbor random walk and those for the (1,2) random walk are indeed the `skipped points' in \cite{w19}.
\begin{definition}
  For (2,1) random walk $X,$ if  $\#\{n\ge0:X_n=k\}=1,$ we call $k$ a cutpoint of $X,$ while for (1,2) random walk $Y,$ if $\#\{n\ge0:X_n=k\}=0,$ we call $k$ a cutpoint of $Y.$
\end{definition}

For $n\ge2,$ let $\varrho_k$ be the one in \eqref{pk} and set
\begin{align}\label{dxyn}D_X(n)=1+\sum_{j=n+1}^{\infty}\prod_{i=n+1}^{j}\varrho_i \text{ and }D_Y(n)=1+\sum_{j=n+1}^{\infty}\prod_{i=n+1}^{j}\varrho_i^{-1}.
\end{align}
  The theorem below gives criteria for the finiteness of the number of cutpoints for both the chain $X$ and $Y.$

\begin{theorem}\label{main}  Consider the chain $Z\in \{X,Y\}.$
Suppose that Condition (C) holds, $a+b=1$ and there exists a number $N_0>0$ such that $\varrho_k$ is increasing in $k\ge N_0$ when $Z=X$ and decreasing in $k\ge N_0$ when $Z=Y.$
  If $$\sum_{n=2}^\infty\frac{1}{D_Z(n)\log n}<\infty,$$ then almost surely, the Markov chain $Z$ has at most finitely many cutpoints.
   If there exists some $\delta>0$ such that $D_Z(n)\le \delta n\log n$ for $n$ large enough and $$\sum_{n=2}^\infty\frac{1}{D_Z(n)\log n}=\infty,$$
 then  almost surely, the Markov chain $Z$ has infinitely many cutpoints.
      \end{theorem}
\begin{remark} \textbf{(i)} For $(1,2) $ random walk $Y,$ we generalize the results of \cite{w19} in the following aspects. On one hand, in \cite{w19}, for the divergent part, it only shown that with a probability $p>2/3,$ the chain $Y$ has infinitely many cutpoints. In our setting, we get an almost-sure result for the divergent case. For this  purpose, asymptotics of the hitting probabilities and the escape probabilities  for (1,2) random walk $Y$ we studied in Proposition \ref{hp0} play a key role.
   On the other hand, in our setting, we allow the existence of size-one right-oriented jumps. But in \cite{w19}, it is required that the right-oriented jumps are always with size 2.  In addition, the criteria are given by $D_Z(n),Z\in\{X,Y\},$ which is defined in terms of $\varrho_k,k\ge2.$ Clearly, $\varrho_k,k\ge2$ are directly computable.
  In \cite{w19}, the criteria rely on the asymptotics of the tail of continued fraction, which is indeed hard to estimate.

  \textbf{(ii)} The proof of Theorem \ref{main} for both the chain $X$ and $Y$ relies essentially on the asymptotics of entries of product of 2-by-2 matrices characterized in terms of products of the tails of a continued fraction, which we studied in another paper \cite{w17}.

\end{remark}
The following corollary gives sharp criteria for finiteness of the cutpoints for the chain $X$ and $Y.$
Fix $\beta\geq0$  and set $$
r_n=\left\{\begin{array}{ll}
  \frac{1}{3}\z(\frac{1}{n}+\frac{1}{n(\log\log n)^\beta}\y),& \text{if } n\ge4,\\
  r_4, & \text{if } n=2,3.
\end{array}\right.
$$
\begin{corollary}\label{tr}
  Assume that Condition (C) holds. \textbf{(i)}
    Suppose $\varrho_k=1-3r_k+\textrm{O}(r_{k}^{2})$ as $k\rto.$ Then
  if $\beta>1,$ $X$ has at most finitely many cutpoints almost surely;
  if $\beta\le 1,$   $X$ has infinitely many cutpoints almost surely. \textbf{(ii)} Suppose $\varrho_k=1+3r_k+\textrm{O}(r_{k}^{2})$ as $k\rto.$ Then
  if $\beta>1,$ $Y$ has at most finitely many cutpoints almost surely;
  if $\beta\le 1,$   $Y$ has infinitely many cutpoints almost surely.
\end{corollary}
\begin{remark}
  The perturbation $r_n,n\ge2$ are constructed based on the observation in \cite{cfrb}. It makes the random walks near-recurrent.
\end{remark}
\proof If $\varrho_k=1\pm3r_k+\textrm{O}(r_{k}^{2})$ as $k\rto,$ then there exists $N_0>0$ such that $\varrho_k$ is decreasing(increasing) in $k\ge N_0.$ Furthermore, with $0<c_1<c_2<\infty$ some proper constants, similar to the proof of  \cite[Theorem 5.1]{cfrb}, if $\varrho_k=1-3r_k+\textrm{O}(r_{k}^{2})$ as $k\rto,$ it can be shown that for $n$ large enough we have \begin{equation}c_1n(\log\log n)^\beta\le D_X(n)\le c_2n(\log\log n)^\beta;\no\end{equation}
otherwise, if $\varrho_k=1+3r_k+\textrm{O}(r_{k}^{2})$ as $k\rto,$ for $n$ large enough we have
\begin{equation}\label{dccy}c_1n(\log\log n)^\beta\le D_Y(n)\le c_2n(\log\log n)^\beta.\end{equation}
Therefore the corollary is a direct consequence of Theorem \ref{main}. \qed

After giving the criteria for the finiteness of the number of cutpoints, another interesting question is to estimate exactly the number of cutpoints in $[0,n]$ when the random walk has infinitely many cutpoints. Intuitively,  the more transient the random walk is, the less  the cutpoints lie in $[0,n].$ We have the following theorem for the near-recurrent random walk.

\begin{theorem}\label{nc} Consider the chain $Z\in \{X,Y\}.$ For $n\ge 2,$ set $S_n=\#\{k\in[2,n]:k\text{ is a cutpoint of }Z\}.$
Suppose that Condition (C) holds. Assume further $\varrho_k=1-3r_k+\textrm{O}(r_{k}^{2})$ as $k\rto$ if $Z=X$ and $\varrho_k=1+3r_k+\textrm{O}(r_{k}^{2})$ as $k\rto$ if $Z=Y$ respectively. If $0\le \beta<1,$ then with proper constants $0<c_3<c_4<\infty,$ we have
\begin{align}\label{lsn}
  c_3\le\liminf_{n\rto} \frac{ ES_n}{\log n(\log\log n)^{-\beta}} \le \limsup_{n\rto} \frac{ ES_n}{\log n(\log\log n)^{-\beta}}\le c_4,
\end{align}
and for each $\varepsilon>0,$
\begin{align}\label{asl}
  \lim_{n\rto} \frac{ S_n}{(\log n)^{1+\varepsilon}(\log\log n)^{-\beta}}=0,\ a.s..
\end{align}
\end{theorem}

\noindent{\bf 1.3 Outline of the paper. } The rest of the paper is organized as follows. Section \ref{pr} gives some preliminary results, including the escape probabilities of the random walks and some basics of continued fractions. In Section \ref{sec3}, by expressing the products of $A_k$'s in terms of the approximants of continued fractions, we study the limit behaviours of the escape probabilities and hitting probabilities of the (1,2) random walk, finishing the proof of Proposition \ref{hp0}. Section \ref{pm} is devoted to proving Theorem \ref{main}. Finally, in Section \ref{pn}, we studying the asymptotics of the number of cutpoints in $[0,n],$ completing the proof of Theorem \ref{nc}.

\section{Preliminary results}\label{pr}
In this section, we compute first the escape probabilities of both the (2,1) and (1,2) random walks, and then introduce some continued fractions and study their tails.

\subsection{Escape probability}\label{sec21}
Let $Q_k^j(m,n)$ and $Q_k(m,n,+)$ be those defined in \eqref{dqa} and \eqref{dqb}. The next lemma gives the escape probabilities of (1,2) random walk.

\begin{lemma}\label{esyp}
 Consider (1,2) random walk $Y$. For $1\le m< k< n,$  we have
\begin{align}\label{epya}&Q_k^n(m,n)=\sum_{s=m+1}^{k}\mathbf e_1 A_s\cdots A_{n-1}\z(\frac{1+\sum_{s=m+1}^{n-1}\mb e_1 A_s\cdots A_{n-1}\mb e_2^t}{1+\sum_{s=m+1}^{n-1}\mb e_1 A_s\cdots A_{n-1}\mb e_1^t}\mb e_1^t-\mb e_2^t\y),\\
  \label{epyb}&Q_k^{n+1}(m,n)=\sum_{s=m+1}^{k}\mathbf e_1 A_s\cdots A_{n-1}\z(\mb e_2^t -\frac{\sum_{s=m+1}^{n-1}\mb e_1 A_s\cdots A_{n-1}\mb e_2^t}{1+\sum_{s=m+1}^{n-1}\mb e_1 A_s\cdots A_{n-1}\mb e_1^t}\mb e_1^t\y),\\
  &Q_k(m,n,+)=\frac{\sum_{s=m+1}^{k}\mb e_1 A_s\cdots A_{n-1}\mb e_1^t}{1+\sum_{s=m+1}^{n-1}\mb e_1 A_s\cdots A_{n-1}\mb e_1^t}\label{epyc}.
\end{align}
\end{lemma}

\proof The proof of  Lemma \ref{esyp} is very standard and is similar to \cite[Lemma 2.1]{br02}. Here, for convenience of the reader, we sketch its proof. Using Markov property, for $1\le m<k<n,j\in \{n,n+1\},$ we have
$
Q_k^j(m,n)=p_{k1}Q_{k+1}^j(m,n)+p_{k2}Q_{k+2}^j(m,n)+q_kQ_{k-1}^j(m,n),
$
which leads to
\begin{align*}
Q_k^j(m,n)-Q_{k-1}^j(m,n)&=\frac{p_{k1}+p_{k2}}{q_k}(Q_{k+1}^j(m,n)-Q_k^j(m,n))\\
&\hspace{0.5cm}+\frac{p_{k2}}{q_k}(Q_{k+2}^j(m,n)-Q_{k+1}^j(m,n)).
\end{align*}
Set
\begin{equation}
V_{k}^{j}=\left(
        \begin{array}{c}
        Q_k^j(m,n)-Q_{k-1}^j(m,n)\\
        Q_{k+1}^j(m,n)-Q_k^j(m,n)  \\
        \end{array}
       \right),\nonumber
\end{equation}
and let $A_k$'s be those defined in (\ref{tn}). Then for $1\le m<k<n,j\in \{n,n+1\},$ we have $V_{k}^{j}=A_kV_{k+1}^j,$ which implies that
$
Q_k^j(m,n)-Q_{k-1}^j(m,n)=\mathbf{e}_{1}A_{k}\cdots A_{n-1}V_{n}^{j}.
$
As a result,
\begin{align}\label{qs}
Q_k^j(m,n)-Q_{m}^j(m,n)=\sum_{s=m+1}^{k}\mathbf{e}_{1}A_{s}\cdot\cdot\cdot A_{n-1}V_{n}^{j}.
\end{align}
Clearly, the boundary values  $Q_{n+1}^n(m,n)=Q_{n}^{n+1}(m,n)=0$ and $ Q_{n+1}^{n+1}(m,n)=Q_{n}^{n}(m,n)=1.$ Thus
\begin{equation}\label{bc}
V_{n}^{n}=\left(
        \begin{array}{c}
        1- Q_{n-1}^n(m,n)\\
        -1  \\
        \end{array}
       \right)
\text{ and }
V_{n}^{n+1}=\left(
        \begin{array}{c}
        - Q_{n-1}^{n+1}(m,n)\\
        1  \\
        \end{array}
       \right).
\end{equation}
Substituting \eqref{bc} into \eqref{qs}, by some careful computation
 we  get \eqref{epya} and \eqref{epyb}.
 In view of  (\ref{dqb}),  taking (\ref{epya}) and (\ref{epyb}) together, we obtain (\ref{epyc}). \qed

Next we compute the escape probabilities of  (2,1) random walk $X$   from certain intervals.  For $0<m-1\le k\le n,$  let
\begin{equation*}
  P_k(m,n,-)=P(X\text{ hits }[0,m]\text{ before it hits }[n,\infty)\big|X_0=k).
\end{equation*}
We have the following lemma.
\begin{lemma}\label{esy}
Consider (2,1) random walk $X$. For $1\le m< k< n,$  we have
\begin{align}\label{epp}
P_k(m,n,-)=\frac{\sum_{s=k}^{n-1}\mathbf e_1 A_s\cdots A_{m+1}\mathbf e_1^t}{1+\sum_{s=m+1}^{n-1}\mathbf e_1 A_s\cdots A_{m+1}\mathbf e_1^t}.\end{align}
\end{lemma}
We remark that \eqref{epp} is just a space reversal argument of \eqref{epyc}, so we omit its proof here.

\subsection{Continued fractions and their tails}\label{sec22}
 We see from Lemma \ref{esyp} and Lemma \ref{esy} that both the escape probabilities of (1,2) and (2,1) random walks are written in terms of products of nonnegative matrices, which are usually hard to compute or estimate. But fortunately,
 product of 2-by-2 matrices is closely related to certain continued fractions so that continued fractions are power tools to  analyze the product of 2-by-2 matrices.

Let $\beta_k,\alpha_k,k\ge 1$ be certain  real numbers. For $1\le k\le n,$ we denote by
\begin{equation}\label{aprx}
\xi_{k,n}\equiv\frac{\beta_k}{\alpha_k}\begin{array}{c}
                                \\
                               +
                             \end{array}\frac{\beta_{k+1}}{ \alpha_{k+1}}\begin{array}{c}
                                \\
                               +\cdots+
                             \end{array}\frac{\beta_n}{\alpha_n}:=\dfrac{\beta_k}{\alpha_k+\dfrac{\beta_{k+1}}{\alpha_{k+1}+_{\ddots_{\textstyle +\frac{\textstyle\beta_{n}}{\textstyle\alpha_{n} } }}}}
\end{equation}
the $(n-k+1)$-th approximant of a continued fraction
 \begin{align}\label{xic}
   &\xi_k:=\frac{\beta_{k}}{\alpha_{k}}\begin{array}{c}
                                \\
                               +
                             \end{array}\frac{\beta_{k+1}}{\alpha_{k+1 }}\begin{array}{c}
                                \\
                               +
                             \end{array}\frac{\beta_{k+2}}{\alpha_{k+2}}\begin{array}{c}
                                \\
                               +\cdots
                             \end{array}.
\end{align}
In the literature,  $\xi_k,k\ge1$ defined in \eqref{xic} are usually called  the tails and $h_k:= \frac{\beta_k}{\alpha_{k}}_{+}\frac{\beta_{k-1}}{ \alpha_{k-1}}_{+\cdots+}\frac{\beta_1}{\alpha_1},$ $k
                             \ge1$
 are usually called the critical tails of the continued fraction ${\frac{\beta_{1}}{\alpha_{1}}}_{+}\frac{\beta_{2}}{\alpha_{2 }}_{+\cdots}$ respectively.
If  $\lim_{n\rto}\xi_{k,n}$ exists, then we say that the continued fraction $\xi_k$  is convergent and its value is defined as $\lim_{n\rto}\xi_{k,n}.$ If
\begin{align}
  \forall k\ge1, \alpha_k>0,\beta_k>0 \text{ and }\exists B>0 \text{ such that } \forall k\ge1,\ B^{-1}\le {\beta_k}/{\alpha_k}\le B,\label{a3}
\end{align}then  by Seidel-Stern theorem(see  \cite[Theorem 3.14]{lw92}), for any $k\ge1,$  $\xi_k$ is convergent. The lemma below will  also be used times and again.

 \begin{lemma}\label{ct}
If $
  \lim_{n\rto}\alpha_n=\alpha\ne0,\ \lim_{n\rto}\beta_n=\beta \text{ and } \alpha^2+4\beta\ge0,$
 then
 for any $k\ge1,$  $\lim_{n\rto}\xi_{k,n}$ exists and furthermore \begin{align}
  \lim_{k\rto}h_k= \lim_{k\rto}\xi_k=\frac{\alpha}{2}\z(\sqrt{1+4\beta/\alpha^2}-1\y).\no
 \end{align}
 \end{lemma}
The proof Lemma \ref{ct} can be found in many references, we refer the reader to \cite{lor}, see  discussion between (4.1) and (4.2) on page 81 therein.

The following lemma gives various inequalities related to the tails and theirs approximants of a continued fraction which are useful for us.

\begin{lemma}\label{xp}
For $1\le k\le n,$ let $\xi_{k,n}$ and $\xi_k$ be those defined in \eqref{aprx} and \eqref{xic}. Suppose that $\alpha_k,\beta_k>0, \forall k\ge1$ and \eqref{a3} is satisfied. Then we have
\begin{align}
 &\xi_{k,n}\rightarrow\xi_k\in(0,\infty), \text{ as } n\rto,\label{xil}\\
  & \xi_{k,n}\Big\{\begin{array}{cc}
     <\xi_k,& \text{if } n-k+1 \text{ is even,} \\
     >\xi_k,&\text{if } n-k+1 \text{ is odd,}
   \end{array},1\le k\le n,\label{pxi}\\
   &\xi_{k,n}\xi_{k+1,n}\Big\{\begin{array}{cc}
     >\xi_k\xi_{k+1},& \text{if } n-k+1 \text{ is even,} \\
     <\xi_k\xi_{k+1},&\text{if } n-k+1 \text{ is odd,}
   \end{array}1\le k\le n-1, \label{txi}\\
  &\xi_k\cdots\xi_n\le \xi_{k,n}\cdots\xi_{n,n} \le \xi_k\cdots\xi_{n-1}\xi_{n,n},1\le k\le n.\label{txii}
  \end{align}
  Furthermore if we assume in addition  $\alpha_k\ge1, \forall k\ge1,$ then
  \begin{align}
    &\xi_{n,n}<\xi_n+\xi_n\xi_{n+1},\label{nxkk}
  \end{align}
  for $1\le k\le n-1,$
 \begin{align}\label{pp}
   &\xi_{k,n}+\xi_{k,n}\xi_{k+1,n}\Big\{\begin{array}{cc}
     \ge\xi_k+\xi_k\xi_{k+1},& \text{if } n-k+1 \text{ is even,} \\
     \le \xi_k+\xi_k\xi_{k+1},&\text{if } n-k+1 \text{ is odd,}
   \end{array}
    \end{align}
    and for $1\le k\le n,$
    \begin{align}
     &\xi_k\cdots\xi_n\le \xi_{k,n}\cdots\xi_{n,n} \le \xi_k\cdots\xi_n(1+\xi_{n+1}),\label{xud} \\ &\sum_{j=k}^{n}\xi_k\cdots\xi_{j}\le\sum_{j=k}^{n}\xi_{k,n}\cdots\xi_{j,n} \le \sum_{j=k}^{n+1}\xi_k\cdots\xi_{j}.\label{sxud}
    \end{align}
  \end{lemma}
\proof Applying Seidel-Stern theorem(see \cite[Theorem 3.14]{lw92}), we get \eqref{xil} and with \eqref{xil} in hand,  \eqref{pxi} is a direct consequence of \cite[Thoerem 3.12]{lw92}. From \eqref{aprx} and \eqref{xic}, we can infer  that
\begin{align*}
  \xi_{k,n}\xi_{k+1,n}=\beta_k-\alpha_k\xi_{k,n} \text{ and } \xi_k\xi_{k+1}=\beta_k-\alpha_k\xi_k.
\end{align*}
As a consequence, we get \begin{align}
  &\xi_{k,n}\xi_{k+1,n}- \xi_k\xi_{k+1}=\alpha_k(\xi_k-\xi_{k,n}),\label{xkk} \\
   &\xi_{k,n}+\xi_{k,n}\xi_{k+1,n}-\xi_k-\xi_k\xi_{k+1}=(\alpha_k-1)(\xi_k-\xi_{k,n})\label{xmk}.
\end{align}
Taking \eqref{pxi} into account, from \eqref{xkk} we get \eqref{txi}, and from \eqref{xmk} we get \eqref{pp} if $\alpha_k\ge 1,\ \forall k\ge1.$
Next, assume $ \alpha_n\ge1,$ $\forall n\ge1.$  Since $\xi_{n+1}>0,$ then
\begin{align*}
  \xi_n&+\xi_n\xi_{n+1}-\xi_{nn}=\xi_n+\beta_n-\alpha_n\xi_{n}-\beta_n/\alpha_n\\
  &=(\alpha_n-1)(\beta_n/\alpha_n-\xi_n)=(\alpha_n-1)(\beta_n/\alpha_n-\beta_n/(\alpha_n+\xi_{n+1}))\ge0,
\end{align*}
which implies \eqref{nxkk}.

For \eqref{txii}, \eqref{xud} and \eqref{sxud}, we prove only the case $k=1$ and $n-k+1$ is an even number. When $k>1$ or $n-k+1$ is odd, the proofs go similarly. Let $n\ge1$ be an even number. Then it follows from \eqref{xil} and \eqref{txi} that
\begin{align*}
  \xi_{1,n}\cdots\xi_{n,n}&=\xi_{1,n} \underline{\xi_{2,n}\xi_{3,n}}\cdots \underline{\xi_{n-2,n}\xi_{n-1,n}}\xi_{n,n}\le \xi_1\underline{\xi_2\xi_3}\cdots\underline{\xi_{n-2}\xi_{n-1}}\xi_{n,n},\\
  \xi_{1,n}\cdots\xi_{n,n}&=\underline{\xi_{1,n}\xi_{2,n}}\cdots \underline{\xi_{n-1,n}\xi_{n,n}}\ge \underline{\xi_1\xi_2}\cdots\underline{\xi_{n-1}\xi_{n}},
\end{align*}
from which we get \eqref{txii}. Putting \eqref{txii} and \eqref{nxkk} together, we obtain \eqref{xud}. To show \eqref{sxud}, on one hand, it follows from \eqref{pxi}-\eqref{pp} that
\begin{align*}
  \sum_{k=1}^{n}&\xi_{1,n}\cdots\xi_{k,n}=\xi_{1,n}+\xi_{1,n}(\xi_{2,n}+\xi_{2,n}\xi_{3,n})+\xi_{1,n}\underline{\xi_{2,n}\xi_{3,n}}(\xi_{4,n}+\xi_{4,n}\xi_{5,n})\\
&+\dots+\xi_{1,n}\underline{\xi_{2,n}\xi_{3,n}}\cdots\underline{\xi_{n-4,n}\xi_{n-3,n}}(\xi_{n-2,n}+\xi_{n-2,n}\xi_{n-1,n})\\
&\quad\quad\quad+\xi_{1,n} \underline{\xi_{2,n}\xi_{3,n}}\cdots \underline{\xi_{n-2,n}\xi_{n-1,n}}\xi_{n,n}\\
&\le \sum_{k=1}^{n+1}\xi_1\cdots\xi_{k},
\end{align*}
which is the upper bound of \eqref{sxud}. On the other hand, by \eqref{txi} and \eqref{pp}, we have
\begin{align*}
  \sum_{k=1}^{n}&\xi_{1,n}\cdots\xi_{k,n}=(\xi_{1,n}+\xi_{1,n}\xi_{2,n})+\underline{\xi_{1,n}\xi_{2,n}}(\xi_{3,n}+\xi_{3,n}\xi_{4,n})\\
&+\dots+\underline{\xi_{1,n}\xi_{2,n}}\cdots\underline{\xi_{n-3,n}\xi_{n-2,n}}(\xi_{n-1,n}+\xi_{n-1,n}\xi_{n,n})\\
&\ge \sum_{k=1}^{n}\xi_1\cdots\xi_{k},
\end{align*} which yields the lower bound of \eqref{sxud}. The lemma is proved.  \qed

\section{Asymptotic behaviours of various escape probabilities and hitting probabilities of the walk}\label{sec3}
\subsection{Product of matrices expressed in terms of the approximants of continued fractions}\label{sec31}

For $2\le k\le n,$ set
\begin{align*}
  & y_{k,n}:=\mathbf e_1 A_k\cdots A_{n}\mathbf e_1^t \text{ and }\zeta_{k,n}:=\frac{y_{k+1,n}}{y_{k,n}};\\
  & z_{k,n}:=\mathbf e_1 A_n\cdots A_k\mathbf e_1^t\text{ and }\theta_{k,n}:=\frac{z_{k+1,n}}{z_{k,n}}.
\end{align*}
 Since the empty product equals identity,  $y_{n+1,n}=z_{n+1,n}=\mathbf e_1I\mathbf e_1^t=1.$ Therefore, we have
\begin{align}
\zeta_{k,n}^{-1}\cdots \zeta_{n,n}^{-1}&=y_{k,n}=\mathbf e_1 A_k\cdots A_{n}\mathbf e_1^t,\label{c4}\\
\theta_{k,n}^{-1}\cdots \theta_{n,n}^{-1}&=z_{k,n}=\mathbf e_1 A_n\cdots A_{k}\mathbf e_1^t,\no
\end{align}
and consequently
\begin{align}
\sum_{k=1}^{n+1} \mathbf e_1A_k\cdots A_n\mathbf e_1^t&=\sum_{k=1}^{n+1} \zeta_{k,n}^{-1}\cdots \zeta_{n,n}^{-1}=\frac{\sum_{k=1}^{n+1} \zeta_{1,n}\cdots \zeta_{k-1,n}}{\zeta_{1,n}\cdots \zeta_{n,n}},\no\\
\sum_{k=1}^{n+1} \mathbf e_1A_n\cdots A_k\mathbf e_1^t&=\sum_{k=1}^{n+1} \theta_{k,n}^{-1}\cdots \theta_{n,n}^{-1}=\frac{\sum_{k=1}^{n+1} \theta_{1,n}\cdots \theta_{k-1,n}}{\theta_{1,n}\cdots \theta_{n,n}}.\no
\end{align}
It is easily seen that
\begin{align}
\theta_{k,n}&=\frac{\mb e_1 A_{n}\cdots A_{k+1}\mb e_1^t}{\mb e_1 A_{n}\cdots A_{k}\mb e_1^t}=\frac{\mb e_1 A_{n}\cdots A_{k+1}\mb e_1^t}{a_k\mb e_1A_{n}\cdots A_{k+1}\mb e_1^t+\mb e_1 A_{n}\cdots A_{k+1}\mb e_2^t}\no\\
&=\frac{\mb e_1 A_{n}\cdots A_{k+1}\mb e_1^t}{a_k\mb e_1A_{n}\cdots A_{k+1}\mb e_1^t+b_{k+1}\mb e_1 A_n\cdots A_{k+2}\mb e_1^t}\no\\
&=\frac{b_{k+1}^{-1}}{a_kb_{k+1}^{-1}+\theta_{k+1,n}}\label{tk}\end{align}
and
\begin{align}
\zeta_{k,n}&=\frac{\mb e_1 A_{k+1}\cdots A_n\mb e_1^t}{\mb e_1 A_{k}\cdots A_n\mb e_1^t}=\frac{\mb e_1 A_{k+1}\cdots A_n\mb e_1^t}{a_k\mb e_1A_{k+1}\cdots A_n\mb e_1^t+b_k\mb e_2 A_{k+1}\cdots A_n\mb e_1^t}\no\\
&=\frac{\mb e_1 A_{k+1}\cdots A_n\mb e_1^t}{a_k\mb e_1A_{k+1}\cdots A_n\mb e_1^t+b_k\mb e_1 A_{k+2}\cdots A_n\mb e_1^t}\no\\
&=\frac{b_{k}^{-1}}{a_kb_{k}^{-1}+\zeta_{k+1,n}}.\label{zk}
\end{align}
Iterating \eqref{tk} and \eqref{zk}, we obtain
\begin{align}
   \theta_{k,n}=\frac{b_{k+1}^{-1}}{a_kb_{k+1}^{-1}}\begin{array}{c}
                                \\
                               +
                             \end{array}\frac{b_{k+2}^{-1}}{ a_{k+1}b_{k+2}^{-1}}\begin{array}{c}
                                \\
                               +\cdots+
                             \end{array}\frac{b_{n+1}^{-1}}{a_nb_{n+1}^{-1}}\no
 \end{align}
 and
 \begin{align}\label{cfz}
   \zeta_{k,n}=\frac{b_k^{-1}}{a_kb_k^{-1}}\begin{array}{c}
                                \\
                               +
                             \end{array}\frac{b_{k+1}^{-1}}{ a_{k+1}b_{k+1}^{-1}}\begin{array}{c}
                                \\
                               +\cdots+
                             \end{array}\frac{b_n^{-1}}{a_nb_n^{-1}}
 \end{align}
 respectively. Under Condition (C), applying Lemma \ref{ct}, we conclude that for each $k\ge2,$ both \begin{align}\label{dtz}
   \theta_k:=\lim_{n\rto}\theta_{k,n}\text{ and }\zeta_k:=\lim_{n\rto}\zeta_{k,n}
 \end{align} exist  and furthermore
 \begin{align}\label{ltz}
   \lim_{k\rto}\theta_k=\varrho^{-1}=\lim_{k\rto}\zeta_k.
 \end{align}
 For $n> m\ge1,$ introduce notations
\begin{align}
F_X(m,n)&:= 1+\sum_{s=m+1}^{n-1}\mathbf{e}_1 A_s\cdots A_{m+1}\mathbf{e}_1^t,  \label{dfmn}\\
 F_Y(m,n)&:= 1+\sum_{j=m+1}^{n-1}\prod_{i=m+1}^{j}{\zeta_i}, \label{dymn}
\end{align}
and let \begin{equation}\label{dmi}
  F_Z(m):=\lim_{n\rightarrow\infty}F_Z(m,n), Z\in \{X,Y\},
\end{equation}
where the limits exist (possibly being $\infty$) because for fixed $m\ge1,$ both $F_X(m,n)$ and $F_Y(m,n)$ are monotone increasing in $n\ge m.$ Then, from Lemma \ref{esy}, we obtain that
 for $1\le m <k < n,$
$$
P_k(m,n,-)=1-\frac {F_X(m,k)}{F_X(m,n)},
$$ and  especially, for $n\geq 1$,
\begin{equation}\label{qfg}
P_{n+1}(n,\infty,-)=1-\frac {1}{F_X(n)}.
\end{equation}
On the other hand, substituting \eqref{c4} into \eqref{epyc}, we see that for $1\le m<n,$
\begin{align}
Q_{m+1}(m,n,+)&=\frac{\zeta_{m+1,n-1}^{-1}\cdots \zeta_{n-1,n-1}^{-1}}{1+\sum_{s=m+1}^{n-1}\zeta_{s,n-1}^{-1}\cdots \zeta_{n-1,n-1}^{-1}}\no\\
&=\frac{1}{1+\sum_{s=m+1}^{n-1}\zeta_{m+1,n-1}\cdots \zeta_{s,n-1}}.\label{c7}
\end{align}
Evidently, $\alpha_i=a_i/b_i\ge1,\forall i\ge1.$ Thus, it follows from \eqref{c7} and \eqref{sxud} that
\begin{align}
 \frac{1}{1+\sum_{s=m+1}^{n}\zeta_{m+1}\cdots \zeta_{s}}\le Q_{m+1}(m,n,+)\le \frac{1}{1+\sum_{s=m+1}^{n-1}\zeta_{m+1}\cdots \zeta_{s}},\no
\end{align}
that is,
\begin{align}
\label{epcb}
 &\frac{1}{F_Y(m,n+1)}\le Q_{m+1}(m,n,+)\le \frac{1}{F_Y(m,n)}.
 \end{align}
Thus, letting $n\rto$ in \eqref{epcb}, we have
 \begin{align}
\label{epcc}
  &Q_{m+1}(m,\infty,+)= \frac{1}{1+\sum_{s=m+1}^{\infty}\zeta_{m+1}\cdots \zeta_{s}}=\frac{1}{F_Y(m)}.
\end{align}
\begin{lemma} \label{varr}
Suppose that Condition (C) holds and $\varrho_k$ is monotone in $k\ge N_0$. Then there are constants $0<c_5<c_6<\infty$ such that for all $m\ge k\geq 2$, $c_5<\frac{\mathbf{e}_1 A_k\cdots A_m\mathbf{e}_1^t}{\varrho_k\cdots\varrho_m}<c_6.$
\end{lemma}
The proof of Lemma \ref{varr} is similar to the one of  \cite[Lemma 2.2]{w18}, so we do not repeat it here.
 Notice that for $n\ge1,$ $F_X(n)$ is written in terms of products of matrices while $F_Y(n)$ is written in terms of products of the tails of a continued fraction. Lemma \ref{varr} enables us to study the asymptotics of $F_X(n)$ and $F_Y(n)$  by those of $D_X(n)$ and $D_Y(n)$ defined in \eqref{dxyn}  respectively. To this end, for $1\le m< n$ set
\begin{align}
  D_Y(m,n)=1+\sum_{j=m+1}^{n-1}\prod_{i=m+1}^{j}\varrho_i^{-1}.\no
\end{align} We have the following lemma.
\begin{lemma}\label{adf}
  Suppose that Condition (C) holds and $\varrho_k$ is monotone in $k\ge N_0$. Then both $D_X(n)$ and $D_Y(n)$ are monotone in $n>N_0$ and  there exist constants $0<c_7<c_8<\infty$ such that for $n>m\ge1,$
  \begin{align}
   & c_7D_X(m)\le F_X(m)\le c_8D_X(m),\label{dyx}\\
   & c_7D_Y(m)\le F_Y(m)\le c_8D_Y(m),\label{dyf}\\
   & c_7D_Y(m,n)\le F_Y(m,n)\le c_8D_Y(m,n). \label{dfyb}
  \end{align}
\end{lemma}
\proof The monotonicity of $D_X(n)$ and $D_Y(n),$ $n\ge N_0$  follows directly from the one of  $\varrho_k, k\ge N_0.$ The inequalities in \eqref{dyx} is a direct consequence of Lemma \ref{varr}. Note that under Condition (C) it follows from \eqref{cfz}-\eqref{ltz} that for some number $0<c<\infty,$ $\zeta_n<c, n\ge2.$ Then, from \eqref{xud}, we get
$$\zeta_{m+1}\cdots\zeta_{j}\le \zeta_{m+1,j}\cdots\zeta_{j,j} \le c\zeta_{m+1}\cdots\zeta_j,\ j\ge m+1.$$
 Thus, on accounting of \eqref{c4}, we have for $n>m\ge1,$
 \begin{align}
   1+c\sum_{j=m+1}^{n-1}\z(\mathbf e_1A_{m+1}\cdots A_j\mathbf e_1^t\y)^{-1}\le  F_Y(m,n)\le1+\sum_{j=m+1}^{n-1}\z(\mathbf e_1A_{m+1}\cdots A_j\mathbf e_1^t\y)^{-1}.\no
 \end{align}
Therefore, using Lemma \ref{varr}, we get \eqref{dfyb}. Letting $n\rto$ in \eqref{dfyb}, we get \eqref{dyf}. \qed

To end this subsection, we state here the following lemma, which will be used when proving Proposition \ref{hp0} and Theorem \ref{main}.
\begin{lemma}\label{uc}
  Suppose that Condition (C) holds. Then,
  \begin{align}
    &\frac{A_{k+1}\cdots A_{k+n}}{\zeta_{k+1}^{-1}\cdots \zeta_{k+n}^{-1}} \rightarrow \frac{\varrho}{\varrho-\sigma}\text{ and }
    \frac{\sum_{s=1}^{n+1}A_{k+s}\cdots A_{k+n}}{\sum_{s=1}^{n+1}\zeta_{k+s}^{-1}\cdots \zeta_{k+n}^{-1}} \rightarrow \frac{\varrho}{\varrho-\sigma}\no
  \end{align}
  uniformly in $k\ge1$ as $n\rto.$
\end{lemma}
 For the proof of Lemma \ref{uc}, we refer the reader to  \cite[Theorem 1]{w17}.

\subsection{Proof of Proposition \ref{hp0}}\label{sec32}
\proof Suppose that Condition (C) holds. According to the Markov property, we have
$\eta_{k,k}(2)=p_{k2}+q_{k}\eta_{k-1,k-1}(1)\eta_{k,k}(2),$
which leads to
\begin{align}
\eta_{k,k}(2)=\frac{b_k}{a_k+\eta_{k-1,k-1}(2)}, k\geq2.\label{a8}
\end{align}
Iterating (\ref{a8}) and using the fact $\eta_{1,1}(2)=0,$ we have for $k\geq2$,
$$\eta_{k,k}(2)=\frac{b_k}{a_{k}}\begin{array}{c}
             \\
             +
             \end{array}\frac{b_{k-1}}{a_{k-1}}\begin{array}{c}
             \\
             +\cdots+
             \end{array}\frac{b_2}{a_2}.$$
Applying Lemma \ref{ct}, we get
$$\lim_{k\rto}\eta_{k,k}(2)=\frac{-a+\sqrt{a^2+4b}}{2}=-\sigma.$$
Next we find the limit of $h_k(2)$ as $k\rto.$ Using again the Markov property, for $k\geq1$, we have
\begin{align}
h_{k+1}(2)=&h_k(2)\eta_{2k+1,2k+1}(2)+h_k(1)\eta_{2k,2k}(1)\eta_{2k+1,2k+1}(2)\nonumber\\
=&h_k(2)\eta_{2k,2k}(2)\eta_{2k+1,2k+1}(2)+\eta_{2k,2k}(1)\eta_{2k+1,2k+1}(2).\label{a11}
\end{align}
Iterating (\ref{a11}) and using the fact $h_{1}(2)=0$, we get
\begin{equation}
h_{k+1}(2)=\sum_{j=1}^{k}\eta_{2j,2j}(1)\eta_{2j+1,2j+1}(2)\cdot\cdot\cdot\eta_{2k+1,2k+1}(2).\label{a12}
\end{equation}
Writing $V_i=\eta_{2i,2i}(1)\eta_{2i+1,2i+1}(2)$ and $U_i=\eta_{2i,2i}(2)\eta_{2i+1,2i+1}(2).$ Then $0\leq U_i\leq1,0\leq V_i\leq1$, and $$\lim_{i\rightarrow\infty}V_i=-\sigma(1+\sigma),\
\lim_{i\rightarrow\infty}U_i=\sigma^2.$$
Thus for arbitrary $0<\varepsilon<(1+\sigma(1+\sigma))\wedge (1-\sigma^2),$ there exists a number $N>0$ such that $V_i<-\sigma(1+\sigma)+\varepsilon$
and $U_i<\sigma^2+\varepsilon$ for $i>N.$
Therefore, by (\ref{a12}) we have for $k>N$,
\begin{align*}
h_{k+1}(2)&=V_k+\sum_{i=1}^{k-1}V_iU_{i+1}\cdots U_k\\
=&V_k+\sum_{i=1}^{N}V_iU_{i+1}\cdots U_k+\sum_{i=N+1}^{k-1}V_iU_{i+1}\cdots U_k\\
\leq &N (\sigma^2+\varepsilon)^{k-N}+(-\sigma(1+\sigma)+\varepsilon)\sum_{i=N+1}^{k}(\sigma^2+\varepsilon)^{k-i}. \end{align*}
Taking the upper limit, we have
$\limsup_{k\rightarrow\infty}h_k(2)\leq\frac{-\sigma(1+\sigma)+\varepsilon}{1-(\sigma^2+\varepsilon)}.$
Since $\varepsilon$ is arbitrary, letting $\varepsilon\rightarrow 0$,  we get
 $$\limsup_{k\rightarrow\infty}h_k(2)\leq-\frac{\sigma}{1-\sigma}.$$
 A similar argument also yields that $$\liminf_{k\rightarrow\infty} h_k(2)\geq -\frac{\sigma}{1-\sigma}.$$
Consequently, we have $\lim_{k\rightarrow\infty}h_k(2)=-\frac{\sigma}{1-\sigma}.$ We thus finish the proof of \eqref{lhe}.

In order to prove \eqref{ql}, we fix $k\ge1$ and set
\begin{equation}
f^{(k)}_{n}:=\frac{\mathbf{e}_{1}A_{k+1}\cdots A_{k+n}\mathbf{e}_{2}^t}{\mathbf{e}_{1}A_{k+1}\cdots A_{k+n}\mathbf{e}_{1}^t} \text{ and } H_{n}^{(k)}:=\sum_{s=k+1}^{k+n}\mb e_1\prod_{i=s}^{k+n}A_{i}(f_{n}^{(k)}\mathbf e_1^t-\mb e_2^t), n\ge 1.\no
\end{equation}
By \eqref{epya}, we have
\begin{align}\label{qhf}
  &Q_{k+1}^{k+n}(k,k+n)=\mathbf e_1 A_{k+1}\cdots A_{k+n-1}\z(\frac{1+\sum_{s=k+1}^{k+n-1}\mb e_1 A_s\cdots A_{k+n-1}\mb e_2^t}{1+\sum_{s=k+1}^{k+n-1}\mb e_1 A_s\cdots A_{k+n-1}\mb e_1^t}\mb e_1^t-\mb e_2^t\y)\no\\
  &\ =\frac{\mathbf e_1 A_{k+1}\cdots A_{k+n-1}\mb e_1^t\z(1-f_{n-1}^{(k)}+\sum_{s=k+1}^{k+n-1}\mb e_1 A_s\cdots A_{k+n-1}(\mb e_2^t-f_{n-1}^{(k)}\mb e_1^t)\y)}{1+\sum_{s=k+1}^{n-1}\mb e_1 A_s\cdots A_{k+n-1}\mb e_1^t}\no\\
  &\ =\frac{\mathbf e_1 A_{k+1}\cdots A_{k+n-1}\mb e_1^t}{1+\sum_{s=k+1}^{k+n-1}\mb e_1 A_s\cdots A_{k+n-1}\mb e_1^t}\z(1-f_{n-1}^{(k)}-H_{n-1}^{(k)}\y).
\end{align}
Note that
\begin{align}
  \label{fr}f_n^{(k)}&=\frac{\mb e_1A_{k+1}\cdots A_{k+n}\mathbf e_2^t}{\mb e_1A_{k+1}\cdots A_{k+n}\mathbf e_1^t}\no\\
  &=\frac{ b_{k+n}\mb e_1A_{k+1}\cdots A_{k+n-1}\mb e_1^t}{\mb e_1A_{k+1}\cdots A_{k+n-1}( a_{k+n}\mathbf e_1^t+\mb e_2^t)}=\frac{ b_{k+n}}{a_{k+n}+f_{n-1}^{(k)}},
\end{align}
which leads to
\begin{align}
  f_n^{(k)}f_{n-1}^{(k)}= b_{k+n}-a_{k+n}f_{n}^{(k)},n\ge 1.\no
\end{align}
Consequently, for $n\ge 1,$
\begin{align}
  H_n^{(k)}&=\sum_{s=k+1}^{k+n}\mb e_1\prod_{i=s}^{k+n}A_{i}(f_{n}^{(k)}\mathbf e_1^t-\mb e_2^t)\label{hr}\\
  &=a_{k+n}f_n^{(k)}- b_{k+n}+\sum_{s=k+1}^{k+n-1}\mb e_1\prod_{i=s}^{k+n-1}A_{i}((a_{k+n}f_n^{(k)}- b_{k+n})\mathbf e_1^t+f_n^{(k)}\mb e_2^t)\no\\
  &=- f_n^{(k)}f_{n-1}^{(k)} -f_n^{(k)}\sum_{s=k+1}^{k+n-1}\mb e_1\prod_{i=s}^{k+n-1}A_{i}(f_{n-1}^{(k)}\mathbf e_1^t-\mb e_2^t)\no\\
  &=- f_n^{(k)}f_{n-1}^{(k)} -f_n^{(k)} H_{n-1}^{(k)}.\no
\end{align}
Since $f_{1}^{(k)}=b_{k+1}/a_{k+1}$ and $H_{1}^{(k)}=0,$ iterating \eqref{hr}, we get
\begin{align}\label{h}
  H_n^{(k)}=\sum_{s=1}^{n-1}(-1)^{n-s}f_s^{(k)} f_{s+1}^{(k)}\cdots  f_n^{(k)},n\ge 1,
\end{align}
and iterating \eqref{fr}, we get
\begin{align}
  f_{n}^{(k)}=\frac{ b_{n+k}}{a_{n+k}}\begin{array}{c}
                                \\
                               +
                             \end{array}\frac{ b_{n-1+k}}{  a_{n-1+k}}\begin{array}{c}
                                \\
                               +\cdots+
                             \end{array}\frac{ b_{k+1}}{ a_{k+1}},n\ge 1.\label{fnra}
\end{align}
We claim that  \begin{align}\label{fu}
  \lim_{n\rto}f_{n}^{(k)}=-\sigma\text{ uniformly in }k.
\end{align} Indeed, note that
\begin{align*}
  \frac{f_{n}^{(k)}}{b_{n+k}}&=\frac{\mb e_1A_{k+1}\cdots A_{k+n-1}\mb e_1^t}{\mb e_1A_{k+1}\cdots A_{k+n}\mb e_1^t}\no\\
 & =\frac{\mb e_1A_{k+1}\cdots A_{k+n-1}\mb e_1^t}{\zeta_{k+1}^{-1}\cdots \zeta_{k+n-1}^{-1}}\frac{\zeta_{k+1}^{-1}\cdots \zeta_{k+n}^{-1}}{\mb e_1A_{k+1}\cdots A_{k+n}\mb e_1^t}\zeta_{k+n}.
\end{align*}
Then it follows from \eqref{ltz} and Lemma \ref{uc} that \eqref{fu} holds.

Substituting \eqref{h} into \eqref{qhf}, we obtain
\begin{align}\label{qht}
  &Q_{k+1}^{k+n}(k,k+n)=\frac{\mathbf e_1 A_{k+1}\cdots A_{k+n-1}\mb e_1^t}{1+\sum_{s=k+1}^{k+n-1}\mb e_1 A_s\cdots A_{k+n-1}\mb e_1^t}\no\\
  &\quad\quad\quad\quad\quad\quad\quad\quad\times\Big(1-f_{n-1}^{(k)} -\sum_{s=1}^{n-2}(-1)^{n-1-s}f_s^{(k)}\cdots  f_{n-1}^{(k)}\Big)\no\\
  &\quad\quad=\frac{\mathbf e_1 A_{k+1}\cdots A_{k+n-1}\mb e_1^t}{1+\sum_{s=k+1}^{k+n-1}\mb e_1 A_s\cdots A_{k+n-1}\mb e_1^t}\Big(1 -\sum_{s=1}^{n-1}(-1)^{n-1-s}f_s^{(k)}\cdots  f_{n-1}^{(k)}\Big).
\end{align}
It is easily seen from \eqref{fnra} that for $n\ge 1, k\ge1,$ $f_n^{(k)}<b_{n+k}/a_{n+k}.$ Since $b_m/a_m\rightarrow b/a$ as $m\rto,$ then for some constant $\kappa>1,$ $\sup_{n\ge2}b_n/a_n<\kappa.$ Thus
$f_n^{(k)}<\kappa$ for all $n\ge1,k\ge1.$
Furthermore, it follows from \eqref{fu} that
 for $0<\varepsilon<1+\sigma,$ $\exists N$ independent of $k$ such that for $s>N,$ $-\sigma-\varepsilon<f_s^{(k)}<-\sigma+\varepsilon.$ With such a  number $N,$ we can write
\begin{align}
  \sum_{s=1}^{n-1}&(-1)^{n-1-s}f_s^{(k)}\cdots  f_{n-1}^{(k)}\no\\
  &=\sum_{s=1}^{N}(-1)^{n-1-s}f_s^{(k)}\cdots  f_{n-1}^{(k)}
  +\sum_{s=N+1}^{n-1}(-1)^{n-1-s}f_s^{(k)}\cdots  f_{n-1}^{(k)}.\label{sf}
\end{align}
But
\begin{align}
   \z|\sum_{s=1}^{N}(-1)^{n-1-s}f_s^{(k)}\cdots  f_{n-1}^{(k)}\y|<\sum_{s=1}^{N}f_s^{(k)}\cdots  f_{n-1}^{(k)}<\kappa^{N}(-\sigma+\varepsilon)^{n-1-N}\rightarrow 0 \label{az}
 \end{align} as $n\rto$ and uniformly in $k\ge1,$
\begin{align*}
\limsup_{n\rto}\sum_{s=N+1}^{n-1}(-1)^{n-1-s}f_s^{(k)}\cdots  f_{n-1}^{(k)}\le \frac{-\sigma+\varepsilon}{1-(-\sigma+\varepsilon)^2}-\frac{(-\sigma-\varepsilon)^2}{1-(-\sigma-\varepsilon)^2}.
\end{align*}
Since $\varepsilon$ is arbitrary, letting $\varepsilon\rightarrow 0,$ we get
\begin{align}
  \limsup_{n\rto}\sum_{s=N+1}^{n-1}(-1)^{n-1-s}f_s^{(k)}\cdots  f_{n-1}^{(k)}\le \frac{-\sigma}{1-\sigma},\label{upf}
\end{align}
uniformly in $k\ge1.$
Taking \eqref{sf}, \eqref{az} and \eqref{upf} together, we obtain \begin{align*}
  \limsup_{n\rto}\sum_{s=1}^{n-1}(-1)^{n-1-s}f_s^{(k)}\cdots  f_{n-1}^{(k)}\le \frac{-\sigma}{1-\sigma},
\end{align*}
uniformly in $k\ge1.$ A similar argument also yields that, uniformly in $k\ge1,$ \begin{align*}
  \liminf_{n\rto}\sum_{s=1}^{n-1}(-1)^{n-1-s}f_s^{(k)}\cdots  f_{n-1}^{(k)}\ge \frac{-\sigma}{1-\sigma}.
\end{align*}
Consequently  we have uniformly in $k\ge1,$
\begin{align}\label{skf}
  \lim_{n\rto}\sum_{s=1}^{n-1}(-1)^{n-1-s}f_s^{(k)}\cdots  f_{n-1}^{(k)}= \frac{-\sigma}{1-\sigma}.
\end{align}
As a result, we can now deduce from \eqref{epya}-\eqref{epyc}, \eqref{qht} and \eqref{skf} that
$$\lim_{n\rto}\frac{Q_{k+1}^{k+n}(k,k+n)}{Q_{k+1}^{k+n+1}(k,k+n)}=-\sigma^{-1}$$ uniformly in $k\ge1,$ which proves  the first limit in \eqref{ql}.

 Finally, we show the second limit in \eqref{ql}.  To this end, note that from \eqref{epyb} we get for $k\ge1,n\ge2,$
\begin{align}
Q_{k+n-1}^{k+n+1}&(k,k+n)\no\\
&=\sum_{s=k+1}^{k+n-1}\mathbf e_1 A_s\cdots A_{k+n-1}\z(\mb e_2^t -\frac{\sum_{s=k+1}^{k+n-1}\mb e_1 A_s\cdots A_{k+n-1}\mb e_2^t}{1+\sum_{s=k+1}^{k+n-1}\mb e_1 A_s\cdots A_{k+n-1}\mb e_1^t}\mb e_1^t\y)\no\\
&=\frac{\sum_{s=k+1}^{k+n-1}\mathbf e_1 A_s\cdots A_{k+n-1}\mb e_2^t}{1+\sum_{s=k+1}^{k+n-1}\mb e_1 A_s\cdots A_{k+n-1}\mb e_1^t}\no\\
&=\frac{b_{k+n-1}\sum_{s=k+1}^{k+n-1}\mathbf e_1 A_s\cdots A_{k+n-2}\mb e_1^t}{1+\sum_{s=k+1}^{k+n-1}\mb e_1 A_s\cdots A_{k+n-2}(a_{k+n-1}\mathbf{e}_1^t+\mathbf{e}_2^t)}\no\\
&=\frac{b_{k+n-1}}{\frac{1}{\sum_{s=k+1}^{k+n-1}\mathbf e_1 A_s\cdots A_{k+n-2}\mb e_1^t}+a_{k+n-1}+\frac{\sum_{s=k+1}^{k+n-2}\mb e_1 A_s\cdots A_{k+n-2}\mathbf{e}_2^t}{1+\sum_{s=k+1}^{k+n-2}\mathbf e_1 A_s\cdots A_{k+n-2}\mb e_1^t}}\no\\
&=\frac{b_{k+n-1}}{\frac{1}{\sum_{s=k+1}^{k+n-1}\mathbf e_1 A_s\cdots A_{k+n-2}\mb e_1^t}+a_{k+n-1}+Q_{k+n-2}^{k+n}(k,k+n-1)}
\label{a16}.
\end{align}
For $m\ge 1,$ set $\beta_{k+m}=b_{k+m}$ and $\alpha_{k+m} =a_{k+m}+\frac{1}{\sum_{s=k+1}^{k+m}\mathbf e_1 A_s\cdots A_{k+m-1}\mb e_1^t}.$ Then iterating \eqref{a16} and using the fact $Q_k^{k+2}(k,k+1)=0,$  we have
\begin{align}\label{qc}
  Q_{k+n-1}^{k+n+1}(k,k+n)=\frac{\beta_{k+n-1}}{\alpha_{k+n-1}}\begin{array}{c}
                                \\
                               +
                             \end{array}\frac{\beta_{k+n-2}}{ \alpha_{k+n-2}}\begin{array}{c}
                                \\
                               +\cdots+
                             \end{array}\frac{\beta_{k+1}}{\alpha_{k+1}},k\ge1,n\ge2.
\end{align}
If we set for $i\ge2,$ $B_i=\left(
              \begin{array}{cc}
                \alpha_i & 1 \\
                \beta_i & 0 \\
              \end{array}
            \right),$ then by \eqref{qc} it is easy to check that
\begin{align}
  Q_{k+n-1}^{k+n+1}(k,k+n)=\frac{\mathbf e_2B_{k+n-1}\cdots B_{k+1}\mathbf e_1^t}{\mathbf e_1B_{k+n-1}\cdots B_{k+1}\mathbf e_1^t}.\no
\end{align}
Notice that by Lemma \ref{uc}, we get
\begin{align}
  &\z(\sum_{s=k+1}^{k+m}\mathbf e_1 A_s\cdots A_{k+m-1}\mb e_1^t\y)^{-1}\sim \frac{\varrho-\sigma}{\varrho} \frac{1}{\sum_{s=k+1}^{k+m}\zeta_{s}^{-1}\cdots\zeta_{k+m-1}^{-1}}\no\\
    &\quad\quad=\frac{\varrho-\sigma}{\varrho}\frac{1}{\zeta_{k+m}} \frac{\zeta_{k+1}\cdots\zeta_{k+m}}{\sum_{s=k+1}^{k+m}\zeta_{k+1}\cdots\zeta_{s-1}}
    \rightarrow\left\{\begin{array}{ll}
                        0, & \text{if }\varrho\ge1, \\
                        \frac{(1-\varrho)(\varrho-\sigma)}{\varrho}, & \text{if }\varrho<1,
                      \end{array}
    \right.\no
\end{align}
uniformly in $k\ge1,$ as $m\rto,$ where  the last convergence is a consequence of \cite[Lemma 4.2]{wy22}.
Thus, we have $\lim_{m\rto}\beta_m=b$ and $\lim_{m\rto}\alpha_m=a+\frac{(1-\varrho)(\varrho-\sigma)}{\varrho}1_{\{\varrho<1\}}.$ Therefore, it follows by Lemma \ref{ct} and \cite[Theorem 1]{w17} that
$$Q_{k+n-1}^{k+n+1}(k,k+n)\rightarrow \tau:=
\left\{\begin{array}{ll}-\sigma, &\text{if }\varrho\ge1,\\
\frac{\sqrt{\hat a^2+4b}-\hat a}{2},& \text{if }\varrho<1,\end{array}\right.$$  uniformly in $k\ge 1$ as $n\rto.$ We thus obtain the second limit in \eqref{ql}. \qed

\subsection{Proof of Proposition \ref{rt}}\label{sec33}
\proof
{\bf(i)} Suppose that $X$ is transient. Then  we  have $P_{m+1}(m,\infty,-)<1$ for any $m\geq1.$ Thus, it follows from (\ref{qfg}) that
$
 F_X(m)=1+\sum_{s=m+1}^\infty\mathbf{e}_1A_s\cdots A_{m+1}\mathbf{e}_1^t<\infty.
 $
 Therefore, from \eqref{dyx}, we get $D_X(m)<\infty$ for all $m\ge1.$
Especially, we have $\sum_{s=2}^{\infty}\varrho_{2}\cdots \varrho_{s}=D_X(1)-1<\infty.$ On the other hand, if  $\sum_{s=2}^{\infty}\varrho_{2}\cdots \varrho_{s}<\infty,$ then  $D(1)=1+\sum_{j=2}^{\infty}\prod_{i=2}^{j}\varrho_i<\infty$. But we see from \eqref{pk} that for all $k\ge2,$ $0<\varrho_k<\infty.$ Thus we have $D_X(m)<\infty$ for all $m\ge1.$  Using \eqref{dyx}, we conclude that
$F_X(m)<\infty$ for all $m\ge1.$ As a result, we see from \eqref{qfg} that $P_{m+1}(m,\infty,-)<1$ for any $m\geq1.$ Consequently, the chain $X$ is transient.

{\bf(ii)}  Suppose that $Y$ is transient. Then $1-Q_{m+1}(m,\infty,+)<1$ for all $m\ge1.$ Thus, by (\ref{epcc}) we get $F_Y(m)=1+\sum_{s=m+1}^{\infty}\zeta_{m+1}\cdots\zeta_s<\infty$. Applying \eqref{dyf}, we have
$D_Y(m)<\infty$ for all $m\ge1.$ Therefore,
 $\sum_{s=2}^{\infty}\varrho_{2}^{-1}\cdots \varrho_{s}^{-1}=D_Y(1)-1<\infty.$ Suppose now $\sum_{s=2}^{\infty}\varrho_{2}^{-1}\cdots \varrho_{s}^{-1}<\infty.$ Then $D_Y(1)=1+\sum_{j=2}^{\infty}\prod_{i=2}^{j}\varrho_i^{-1}<\infty.$ Since $0<\varrho_k<\infty$ for all $k\ge2,$ we have $D_Y(m)<\infty$ for all $m\ge1.$ Applying again \eqref{dyf}, we get $F_Y(m)<\infty$ for all $m\ge1.$
Therefore, it follows from (\ref{epcc}) that $1-Q_{m+1}(m,\infty,+)<1$ for all $m\geq1.$ Consequently, $Y$ is transient. \qed

\section{Criteria for finiteness of the number of cutpoints}\label{pm}
The main task of this section is to prove Theorem \ref{main}. It is divided into two subsections. Subsection \ref{tw} and Subsection \ref{wt} are devoted to studying (2,1) random walk and (1,2) random walk respectively.

\subsection{(2,1) random walk-Proof of Theorem \ref{main} for the chain $X$}\label{tw}

Consider the chain $X,$  we divide the proof of Theorem \ref{main} into  two parts, that is,  the convergent part  and the divergent part. To begin with, we compute the probabilities that a site $k$ or both two sites $k$ and $j$ are cutpoints. Let $$ C:=\{k\geq 0:k \text{ is a cutpoint }\}$$ be the collection of all cutpoints of the chain $X.$  For $0\le m<k,$ denote $$C_{m,k}=\{x:2^m<x\leq2^k,x\in C\}$$
and set $A_{m,k}:=\#{C_{m,k}}$. We have the following lemma whose proof is similar to \cite[Lemma 2.1]{cfrb}

\begin{lemma}\label{lem5}
 We have
$$
 P(k\in C)=\frac{q_k}{F_X(k)},k\ge2,
$$  and
$$
 P(j\in C, k\in C)=\frac{q_jq_k}{F_X(j,k)F_X(k)}, 2\le j<k.
$$
\end{lemma}

\subsubsection{Convergent part}

Now we begin to prove the convergent part of Theorem \ref{main} for the chain $X.$

\proof The idea of the proof is taken from \cite{jlp}. Suppose that Condition (C) holds, $a+b=1,$ $\varrho_k$ is increasing in $k\ge N_0$ and $\sum_{n=2}^\infty\frac{1}{D_X(n)\log n}<\infty.$ Then we see from \eqref{ro} that
\begin{align}
\varrho_k \rightarrow 1, \text{ as }k\rto. \label{ru1}
\end{align}
For $m\ge1,$ let $l_m$ be the largest $k\in C_{m,m+1}$ if $C_{m,m+1}$ is nonempty and
write
$$  \chi_m:=\min_{k\in (2^{m},2^{m+1}]}\sum_{i=1}^{2^{m-1}} \frac{1}{F_X(k-i,k)}.
$$
Notice that for $m\ge1,$ $2^m<k\le  2^{m+1}$ and $2^{m-1}< i<k$ we have
\begin{align}
  P(i\in C&|A_{m,m+1}>0,l_m=k)=\frac{P(i\in C,A_{m,m+1}>0,l_m=k)}{P(A_{m,m+1}>0,l_m=k)}\no\\
  &=\frac{P(i\in C,k\in C, j\in C^c,k+1\le j\le 2^{m+1})}{P(k\in C, j\in C^c,k+1\le j\le 2^{m+1})}\no\\
  &=\frac{\frac{q_iq_k}{F_X(i,k)F_X(k)}P(j\in C^c,k+1\le j\le 2^{m+1}|X_0=k+1)}{\frac{q_k}{F_X(k)}P(j\in C^c,k+1\le j\le 2^{m+1}|X_0=k+1)}\no\\
  &=\frac{q_i}{F_X(i,k)}.\no
\end{align}
Moreover, under Condition (C), we have $c<q_i<1$ for all $i\ge2.$
Thus, for  $m\ge1,$  we get
       \begin{align}
        \sum_{j=2^{m-1}+1}^{2^{m+1}}&P(j\in C)=E(A_{m-1,m+1})\no\\
              &\ge \sum_{k=2^{m}+1}^{2^{m+1}}P(A_{m,m+1}>0,l_m=k)\times E(A_{m-1,m+1}|A_{m,m+1}>0,l_m=k)\no\\
              &= \sum_{k=2^{m}+1}^{2^{m+1}}P(A_{m,m+1}>0,l_m=k) \sum_{i={2^{m-1}+1}}^kP(i\in C|A_{m,m+1}>0,l_m=k)\no\\
              &= \sum_{k=2^{m}+1}^{2^{m+1}}P(A_{m,m+1}>0,l_m=k)\z(1+ \sum_{i={2^{m-1}+1}}^{k-1}\frac{q_i}{F_X(i,k)}\y)\no\\
              &\ge cP(A_{m,m+1}>0)\min_{2^{m}< k \le 2^{m+1}} \sum_{i={2^{m-1}+1}}^{k-1} \frac{1}{F_X(i,k)}\no\\
               &\ge cP(A_{m,m+1}>0)\min_{2^{m}< k \le 2^{m+1}} \sum_{i=1}^{2^{m-1}}\frac{1}{F_X(k-i,k)}\no\\
              &=c \chi_m P(A_{m,m+1}>0).\label{pc}
      \end{align}
Since $\varrho_k$ is increasing in $k\ge N_0,$  by \eqref{ru1} and Lemma \ref{varr}, we obtain
\begin{align}
F_X(m,n)&=1+\sum_{s=m+1}^{n-1} \mathbf{e}_1 A_{s}\cdots A_{m+1}\mathbf e_1^t\nonumber\\
&\le c \z(1+\sum_{s=m+1}^{n-1} \varrho_s \cdots \varrho_{m+1}\y)\leq c(n-m),\nonumber
\end{align}
which implies that \begin{align}
  \chi_m=\min_{k\in (2^{m},2^{m+1}]}\sum_{i=1}^{2^{m-1}} \frac{1}{F_X(k-i,k)}\ge c\sum_{i=1}^{2^{m-1}}\frac{1}{i}\ge cm.\label{cm}
\end{align}
Therefore, taking \eqref{dyx} and Lemma \ref{lem5} into consideration, from \eqref{pc} and \eqref{cm}, we get
\begin{align*}
\sum_{m=1}^\infty &P(A_{m,m+1}>0)\le c\sum_{m=1}^\infty\frac{1}{\chi_m}\sum_{j=2^{m-1}+1}^{2^{m+1}}P(j\in C)\\
&\le c\sum_{m=1}^\infty\frac{1}{m}\sum_{j=2^{m-1}+1}^{2^{m+1}} \frac{1}{F_X(j)}
= c\sum_{m=1}^\infty\frac{1}{m}\sum_{j=2^{m-1}+1}^{2^{m+1}} \frac{1}{D_X(j)}\\
&\le c\sum_{m=1}^\infty \sum_{j=2^{m-1}+1}^{2^{m+1}} \frac{1}{D_X(j)\log j}\\
&\le c\sum_{m=1}^\infty \frac{1}{D_X(m)\log m}<\infty.
\end{align*}
An application of Borel-Cantelli lemma yields that with probability one, only finitely many of the events $\{A_{m,m+1}>0\},m\ge1$ occur. We conclude that the  chain $X$ has at most finitely many cutpoints almost surely. The convergent part of Theorem \ref{main} is proved for the chain $X.$ \qed
\subsubsection{Divergent part}
Next, we prove the divergent part of Theorem \ref{main} for the chain $X.$

\proof Suppose that Condition (C) holds, $\varrho_k$ is increasing in $k\ge N_0,$ and moreover, $\exists \delta>0$ such that $D_X(n)\le \delta n\log n$ for $n$ large enough and $\sum_{n=2}^\infty\frac{1}{D_X(n)\log n}=\infty.$ Let $\theta_k, k\ge2$ be those in \eqref{dtz}. For $n> m\ge1,$ set
 \begin{align*}
   G(m,n)&=1+\sum_{s=m+1}^{n-1} \theta_{s}^{-1}\cdots \theta^{-1}_{m+1}=1+\sum_{j=1}^{n-m-1} \theta_{m+j}^{-1}\cdots \theta^{-1}_{m+1},\\
   G(m)&=1+\sum_{s=m+1}^{\infty} \theta_{s}^{-1}\cdots \theta^{-1}_{m+1}=1+\sum_{j=1}^{\infty} \theta_{m+j}^{-1}\cdots \theta^{-1}_{m+1},
 \end{align*}
and let $F_X(m,n)$  be those defined in \eqref{dfmn}.
The following  proposition plays a key role in the proof.
\begin{proposition}\label{dfx} Suppose that Condition (C) holds. Then $\exists N_1>0$ such that for $m\ge N_1$ we have
$\frac{F_X(m,m+n)}{G(m,m+n)}\rightarrow\frac{\varrho}{\varrho-\sigma}$
uniformly in $m$  as $n\rto.$
\end{proposition}
For the proof of Proposition \ref{dfx}, we refer the reader to \cite{w17}. For $k\ge2,$ set \begin{align}
  m_k=[k\log k] \text{ and }E_k=\{m_k\in C\}.
  \no
\end{align} Here and throughout, $[x]$ denotes the integer part of $x.$
To prove the divergent part, we need to show that
\begin{align}
  P(E_k,k\geq 1 \text{ occur infinitely often})=1. \label{pbi}
\end{align}
Fix $\varepsilon>0.$
Applying Proposition \ref{dfx}, we can infer that there exist  $N_1>0$ and $j_1>0$ (independent of $m$) such that
\begin{align}\label{dfb}
   \frac{\varrho}{\varrho-\sigma}-\varepsilon\le \frac{F_X(m,m+n)}{G(m,m+n)}\le \frac{\varrho}{\varrho-\sigma}+\varepsilon, \forall m\ge N_1, n\ge j_1.
 \end{align}
Letting $n\rto$ in \eqref{dfb}, we get
  \begin{align}\label{df}
   \frac{\varrho}{\varrho-\sigma}-\varepsilon\le \frac{F_X(m)}{G(m)}\le \frac{\varrho}{\varrho-\sigma}+\varepsilon, \forall m\ge N_1 .
 \end{align}
Taking \eqref{dfb} and \eqref{df} together, we come to a conclusion that
\begin{align}\label{fg}
\frac{F_X(m,n)}{F_X(m)}\ge \frac{( \frac{\varrho}{\varrho-\sigma}-\varepsilon)G(m,n)}{( \frac{\varrho}{\varrho-\sigma}+\varepsilon)G(m)}, \forall m\ge N_1, n-m\ge j_1.
\end{align}
Also, from \eqref{dyf} and \eqref{df}, we obtain
\begin{align}\label{gd}
  \frac{1}{G(m)}\ge \frac{1}{\varrho/(\varrho-\sigma)-\varepsilon} \frac{1}{F_X(m)}\ge  \frac{1}{\varrho/(\varrho-\sigma)-\varepsilon} \frac{c_8^{-1}}{D_X(m)}=\frac{c_9}{D_X(m)}, m\ge N_1,
\end{align}
where $c_9={(\varrho/(\varrho-\sigma)-\varepsilon)^{-1}c_8^{-1}}.$

Noticing that
 for $l>k,$
$m_l-m_k\ge l\log l-k\log k\ge \log k.$ Thus, there exists a number $k_1>0$ such that for all $k\ge k_1,$
$m_k\ge N_1, m_l-m_k>j_1.$

In order to prove \eqref{pbi}, on one hand,
it follows from \eqref{gd} and Lemma \ref{lem5} that
\begin{align}
  \sum_{k=k_1}^\infty P(E_k)=\sum_{k=k_1}^\infty\frac{q_{m_k}}{F_X(m_k)}\ge \sum_{k=k_1}^\infty\frac{c}{D_X(m_k)}=\sum_{k=k_1}^\infty\frac{c}{D([k\log k])}=\infty,\no
\end{align}
where for the last step, we use \cite[Lemma 2.2]{cfrb}.
On the other hand,  from \eqref{fg} we get
\begin{align}
\frac{F_X(m_k,m_l)}{F_X(m_k)}&\ge \frac{\z(\frac{\varrho}{\varrho-\sigma}-\varepsilon\y)G(m_k,m_l)}{\z(\frac{\varrho}{\varrho-\sigma}+\varepsilon\y)G(m_k)}=\frac{\frac{\varrho}{\varrho-\sigma}-\varepsilon}{\frac{\varrho}{\varrho-\sigma}+\varepsilon}\z(1-\prod_{i=m_k}^{m_l-1}\z(1-\frac{1}{G(i)}\y)\y)\no\\
&\geq\frac{\frac{\varrho}{\varrho-\sigma}-\varepsilon}{\frac{\varrho}{\varrho-\sigma}+\varepsilon}\z(1-\exp\z\{-\sum_{i={m_k}}^{m_l-1}\frac{1}{G(i)}\y\}\y), \forall k\ge k_1. \label{ff}
\end{align}
Let
\begin{equation}\label{btt}
\ell=\min\z\{l\ge k:\sum_{i=m_k}^{m_l-1}\frac{1}{G(i)}\ge \log\frac{1+\varepsilon}{\varepsilon}\y\}.
\end{equation}
Then for $l\geq \ell$,  we have
\begin{align*}
\frac{F_X(m_k)}{F_X(m_k,m_l)}\leq\frac{\frac{\varrho}{\varrho-\sigma}+
\varepsilon}{\frac{\varrho}{\varrho-\sigma}-\varepsilon}
\z(1-\exp\z\{-\sum_{i=m_k}^{m_l-1}\frac{1}{G(i)}\y\}\y)^{-1}\leq
\frac{\frac{\varrho}{\varrho-\sigma}+\varepsilon}{\frac{\varrho}{\varrho-\sigma}
-\varepsilon}(1+\varepsilon).
\end{align*}
Therefore, using Lemma \ref{lem5}, we get
\begin{align}\label{akal1}
P(E_kE_l)&=P(m_k\in C, m_l\in C)=\frac{F_X(m_k)}{F_X(m_k,m_l)}P(E_k)P(E_l)\nonumber\\
&\le \frac{\frac{\varrho}{\varrho-\sigma}+\varepsilon}{\frac{\varrho}{\varrho-\sigma}-\varepsilon}
(1+\varepsilon)P(E_k)P(E_l), \ l\ge \ell.
\end{align}
Next we consider
 $k<l<\ell$. Note that
$
1-e^{-u}\ge cu,$ for $ u\in \z[0,\log\frac{1+\varepsilon}{\varepsilon}\y].
$
Then using \eqref{df}, \eqref{ff}, Lemma \ref{adf} and Lemma \ref{lem5}, we have
\begin{align*}
P(E_kE_l) &= \frac{F_X(m_k)}{F_X(m_k,m_l)}P(E_k)P(E_l)\le \frac{c P(E_k)P(E_l)}{\sum_{i=m_k}^{m_l-1}\frac{1}{G(i)}}\no\\
&\le
 \frac{c P(E_k)P(E_l)}{\sum_{i=m_k}^{m_l-1}\frac{1}{F_X(i)}}\le
 \frac{c P(E_k)P(E_l)}{\sum_{i=m_k}^{m_l-1}\frac{1}{D_X(i)}}.\no\\
  \end{align*}
 Since $D_X(n)$ is increasing in $N\ge N_0,$  thus applying again Lemma \ref{lem5}, we get
 \begin{align*}P(E_kE_l)\le \frac{c P(E_k)P(E_l)D_X(m_l)}{m_l-m_k} \le \frac{c P(E_k)P(E_l)F_X(m_l)}{m_l-m_k}\le \frac{c P(E_k)}{l\log l-k\log k}.
\end{align*}
As a result,
\begin{align}\label{akal2}
\sum_{l=k+1}^{\ell-1}&P(E_kE_l)\le cP(E_k) \sum_{l=k+1}^{\ell-1}\frac{1}{l\log l-k\log k}\no\\
&\le cP(E_k)\frac{1}{\log k}\sum_{l=k+1}^{\ell-1}\frac{1}{l-k}\le  cP(E_k)\frac{\log \ell }{\log k}.
\end{align}
We claim that there exists some $\gamma$ depending only on $\varepsilon$ such that \begin{equation}\label{klg}
  \frac{\log \ell }{\log k}\le\gamma.
\end{equation}
Indeed, by (\ref{btt}), we know that
        \begin{equation*}
          \sum_{i=m_k+1}^{m_l-1}\frac{1}{G(i)}< \log\frac{1+\varepsilon}{\varepsilon}\text{ for }k<l<\ell.
        \end{equation*}
This implies that, for large $k,$ we have $l<k^\gamma$ with $\gamma>\z(\frac{1+\varepsilon}{\varepsilon}\y)^{\delta/c_9}+\varepsilon$ where $c_9$ is the one in \eqref{gd}.  Assume the contrary that $l\ge k^\gamma.$ Recall that by assumption, $D_X(n)<\delta n\log n$ for large $n.$ Then taking \eqref{gd} into account, we get
\begin{align*}
   \sum_{i=m_k+1}^{m_l-1}\frac{1}{G(i)} &\ge c_9 \sum_{i=m_k+1}^{m_l-1}\frac{1}{D_X(i)}\ge \frac{c_9}{\delta}\sum_{i=m_k+1}^{m_l-1}\frac{1}{i\log i}\no\\
   &\ge \frac{c_9}{\delta} \z[\log\log (m_l)-\log\log(m_k+1)\y]\\
   &\ge\frac{c_9}{\delta} \log(\gamma-\varepsilon)\ge\log\frac{1+\varepsilon}{\varepsilon},
\end{align*}
a contradiction. We thus finish the proof of \eqref{klg}.

Now we can deduce from \eqref{akal2} and \eqref{klg} that
\begin{align}
  \sum_{l=k+1}^{\ell-1}P(E_kE_l)\le cP(E_k). \label{sl}
\end{align}
Taking
 (\ref{akal1}) and (\ref{sl}) together, we get
\begin{align*}
\sum_{k=k_1}^{N}\sum_{l=k+1}^NP(E_kE_l) \le\sum_{k=k_1}^{N}\sum_{l=k+1}^N \frac{\frac{\varrho}{\varrho-\sigma}+\varepsilon}{\frac{\varrho}{\varrho-\sigma}-\varepsilon}(1+\varepsilon) P(E_k)P(E_l)+ c\sum_{k=k_1}^{N}P(E_k).
\end{align*}
Writing $H(\varepsilon)=\frac{\frac{\varrho}{\varrho-\sigma}+\varepsilon}{\frac{\varrho}{\varrho-\sigma}-\varepsilon}(1+\varepsilon),$
 we have
\begin{align*}
  \alpha_H&:=\liminf_{N\rto}\frac{\sum_{k=k_1}^{N}\sum_{l=k+1}^NP(E_kE_l)-
  \sum_{k=k_1}^{N}\sum_{l=k+1}^N H P(E_k)P(E_l)}{\z[\sum_{k=k_1}^{N}P(E_k)\y]^2}\\
  &\le \lim_{N\rto}\frac{c}{\sum_{k=k_1}^{N}P(E_k)}=0.
\end{align*}
By a general version of Borel-Cantelli lemma \cite[p. 235]{pe04}, we get
\begin{align*}
P(E_k,k\ge 1 &\text{ occur infinitely often} )
\ge P(E_k,k\ge k_1 \text{ occur infinitely often} )\\
&\ge \frac{1}{H+2\alpha_H}\ge \frac{\frac{\varrho}{\varrho-\sigma}-\varepsilon}{\frac{\varrho}{\varrho-\sigma}+\varepsilon}\frac{1}{1+\varepsilon}.
\end{align*}
Since $\varepsilon$ is arbitrary, letting $\varepsilon\rightarrow 0,$ we conclude that $$P(E_k,k\ge1 \text{ occur infinitely often} )=1.$$
 The divergent part of Theorem \ref{main} for the chain $X$ is proved. \qed

\subsection{(1,2) random walk-Proof of Theorem
\ref{main} for the chain $Y$}\label{wt}

 The basic conception to prove  Theorem
\ref{main} for the chain $Y$ is to split $\mathbb Z_+$ into nonintersecting layers, that is, to let $\mathbb Z^{+}=\bigcup_{k=0}^{\infty}L_{k}$ with $L_{k}=\{2k,2k+1\},k=0,1,2,\cdots$ and then compute the probabilities that $L_{k}$ or both $L_{k}$ and $L_{m}$ contain a cutpoint.

For the chain $Y,$  let $S=\{k\ge 2: k\text{ is a cutpoint}\},$ let
$C^S=\{k\ge1: L_k \text{ contains a cutpoint} \},$ and for $0\le m< k$ denote $C_{m,k}=\{x:2^m<x\leq2^k,x\in C^S\}$
and set $A_{m,k}:=\#{C_{m,k}}.$ With those notations, the following proposition, whose proof will be delayed to the end of this section, is crucial for our proof.
\begin{proposition}\label{p31}
Suppose that  Condition (C) holds and $a+b=1.$ Then we have
\begin{align}
\lim_{k\rightarrow\infty}F_Y(2k)P(k\in C^{S})=-\frac{2\sigma}{1-\sigma},\label{d1}
\end{align}
and for any $\varepsilon>0$, there exist $N>0$ and $k_{0}>0$ such that for $k-N>m>k_{0}$,
\begin{align}
P(m\in C^{S},k\in C^{S})\leq(1+\varepsilon)P(m\in C^{S})P(k\in C^{S})\frac{F_Y(2m+1)}{F_Y(2m+1,2k)}.\label{d2}
\end{align}
\end{proposition}
\begin{remark}
  Compared with \cite[Proposition 3.1]{w19}, the essential improvement we make here is a refinement of the upper bound in \eqref{d2}. In \cite[Proposition 3.1]{w19}, the term corresponding $(1+\varepsilon)$ here is $(3/2+\varepsilon)$ so that in \cite{w19}, for the divergent part, it is only shown that with a probability $p\ge 2/3,$ the chain $Y$ has infinitely many cutpoint, whereas we can get an almost-sure result for the divergent part.
\end{remark}
We also divide the proof of Theorem \ref{main} for the chain $Y$ into two parts. The convergent part is treated in Section \ref{cy}, while the divergent part is proved in Section \ref{doy}.
With Proposition \ref{p31} in hand, the proof of Theorem \ref{main} for $Y$ overlaps the one in \cite{w19} to some extent. For those overlapped contents, we only sketch their proofs.

\subsubsection{Convergent part}\label{cy}
We now prove the convergent part  of Theorem \ref{main} for the chain $Y.$

\proof
Assume that Condition (C) holds, $a+b=1,$ $\varrho_k$ is decreasing in $k\ge N_0$ and furthermore $\sum_{m=2}^\infty \frac{1}{D_Y(n)\log n}<\infty.$
Let $F_Y(m,n)$ and $F_Y(m)$ be those defined in \eqref{dymn} and \eqref{dmi}.  Exactly by the same line of the proof of \cite[Lemma 4.1]{w19}, we have that for $m$ large enough and $k\in C_{m,m+1}, $ $ 2^{m-1}<i\le k,$
     \begin{align}
       &P(i\in C^S|l_m=k,2k\in S)\ge \frac{c}{F_Y(2i,2k+1)},\label{ki1}\\
        & P(i\in C^S|l_m=k,2k+1\in S)\ge \frac{c}{F_Y(2i,2k+1)}.     \label{ki2}
     \end{align}
 Write for $m\ge1,$ $$ \psi_m:=\sum_{i=1}^{2^{m-1}}\min_{2^{m}< k \le 2^{m+1}} \frac{1}{F_Y(2(k-i),2k+1)}.$$
Owing to \eqref{dfyb}, we have
$$\psi_m \ge c\sum_{i=1}^{2^{m-1}}\min_{2^{m}< k \le 2^{m+1}} \frac{1}{D_Y(2(k-i),2k+1)}.$$ Since $\varrho_k$ is decreasing in $k\ge N_0,$ and $\varrho_k^{-1}\rightarrow 1$ as $k\rto,$  then for $1\le m<n,$ $D_Y(m,n)<c(n-m).$ Therefore, we have
$$\psi_m \ge c\sum_{i=1}^{2^{m-1}}\min_{2^{m}< k \le 2^{m+1}} \frac{1}{D_Y(2(k-i),2k+1)}\ge  c\sum_{i=1}^{2^{m-1}}\frac{1}{2i+1}\ge cm,m\ge1. $$
Using \eqref{ki1} and \eqref{ki2}, we have for $m$  large enough,
       \begin{align*}
        \sum_{j=2^{m-1}+1}^{2^{m+1}}P(j\in C^S)
              &\ge cP(A_{m,m+1}>0)\min_{2^{m}< k \le 2^{m+1}} \sum_{i=1}^{2^{m-1}}\frac{1}{F_Y(2(k-i),2k+1)}\no\\
              &=cP(A_{m,m+1}>0)\psi_m\ge cm P(A_{m,m+1}>0).
      \end{align*}
As a result, on accounting of \eqref{dyf} and \eqref{d1}, we get
 \begin{align*}
        \sum_{m=1}^\infty P(A_{m,m+1}>0)&\le\sum_{m=1}^\infty\frac{c}{m}\sum_{j=2^{m-1}+1}^{2^{m+1}}P(j\in C^S) \le \sum_{m=1}^\infty\frac{c}{m}\sum_{j=2^{m-1}+1}^{2^{m+1}} \frac{1}{F_Y(2j)}\no\\
        &\le \sum_{m=1}^\infty\frac{c}{m}\sum_{j=2^{m-1}+1}^{2^{m+1}} \frac{1}{D_Y(2j)}\le \sum_{m=1}^\infty\frac{c}{m}\sum_{j=2^{m}+2}^{2^{m+2}} \frac{1}{D_Y(j)}\no\\
        &\le c\sum_{m=1}^\infty \sum_{j=2^{m}+2}^{2^{m+2}} \frac{1}{D_Y(j)\log j} \le c\sum_{m=1}^\infty \frac{1}{D_Y(n)\log n}<\infty.
      \end{align*}
      Applying Borel-Cantelli lemma, we conclude that  with probability one, only finitely many of the events $\{A_{m,m+1}>0\},m\ge1$ occur. The convergent part of Theorem \ref{main} for the chain $Y$ is proved. \qed

\subsubsection{Divergent part}\label{doy}
We now turn to prove the divergent part of Theorem \ref{main} for the chain $Y.$

\proof
With Proposition \ref{p31} in hand, the proof of the divergent part is more or less similar to the one for the chain $X.$ So we  only sketch its proof here. Suppose that Condition (C) holds, $\varrho_k$ is decreasing in $k\ge N_0,$  for large $n,$ $D_Y(n)<\delta n\log n$ for some $\delta>0$ and $\sum_{n=2}^\infty \frac{1}{D_Y(n)\log n}=\infty.$ Since there is no danger of making confusion, for the chain $Y,$ for $k\ge 1$ we also set
$$m_k=[k\log k], E_k=\{k\in C^S\}.$$
It suffices to show
\begin{align}\label{sfy}
  P(E_k,k\geq 1 \text{ occur infinitely often})=1.
\end{align}
Fix $\varepsilon >0.$ Taking \eqref{dfyb} into account, it follows by Proposition \ref{p31} that there exists $k_2>0$ such that for $k\geq k_2$,
\begin{align}
  \sum_{k=k_2}^{\infty}P(E_k)&\geq \sum_{k=k_2}^{\infty}\frac{c}{F_Y(2m_k+1)}\ge \sum_{k=k_2}^{\infty}\frac{c}{D_Y(2m_k+1)}\no\\
  &\ge \sum_{k=k_2}^{\infty}\frac{c}{D_Y(2[k\log k]+1)}\geq \sum_{k=k_2}^{\infty}\frac{c}{D_Y([2k\log 2k])}=\infty, \label{syi}
\end{align}
where for the last step, we use \cite[Lemma 2.2]{cfrb}.

On the other hand, applying again Proposition \ref{p31}, there exists $k_3\ge k_2$ such that
for $l>k>k_3,$
\begin{align}\label{e3}
P(E_kE_l)
&\leq (1+\varepsilon)P(E_k)P(E_l)\frac{F_Y(2m_k+1)}{F_Y(2m_k+1,2m_l)}\nonumber\\
&=(1+\varepsilon)\z\{\frac{F_Y(2m_k+1,2m_l)}{F_Y(2m_k+1)}\y\}^{-1}P(E_k)P(E_l)\no\\
&=(1+\varepsilon)\z\{1-\prod_{i=2m_k+1}^{2m_l-1}\z(1-\frac{1}{F_Y(i)}\y)\y\}^{-1}P(E_k)P(E_l)\no\\
&\le(1+\varepsilon)\z(1-\exp\z\{-\sum_{i=2m_k+1}^{2m_l-1}\frac{1}{F_Y(i)}\y\}\y)^{-1}P(E_k)P(E_l).
\end{align}
Define
\begin{align}
\ell_1=\min\z\{l\geq k:\sum_{i=2m_k+1}^{2m_l-1}\frac{1}{F_Y(i)}\geq \log \frac{1+\varepsilon}{\varepsilon}\y\}.\no
\end{align}
Then for $l\geq \ell_1$ we have
$\z(1-\exp\z\{-\sum_{i=2m_k+1}^{2m_l-1}\frac{1}{F_Y(i)}\y\}\y)^{-1}\leq 1+\varepsilon.$
Therefore, by (\ref{e3}) we have
\begin{align}
P(E_kE_l)\leq (1+\varepsilon)^2P(E_k)P(E_l), \text{  for } l\geq \ell_1.\label{e6}
\end{align}
Next we consider $k<l<\ell_1.$ Note that for $0\leq u\leq \log \frac{1+\varepsilon}{\varepsilon}$, we have $1-e^{-u}\geq cu$. Since $\varrho_k^{-1}$ is increasing in $k> N_0$ and so is $D_Y(k).$   Thus on accounting of \eqref{dfyb}, (\ref{d1}) and (\ref{e3}), we have
\begin{align}
P(E_kE_l)&\le(1+\varepsilon)\z(1-\exp\z\{-\sum_{i=2m_k+1}^{2m_l-1}\frac{1}{F_Y(i)}\y\}\y)^{-1}P(E_k)P(E_l)\nonumber\\
&\leq\frac{cP(E_k)P(E_l)}{\sum_{i=2m_k+1}^{2m_l-1}\frac{1}{F_Y(i)}}\leq\frac{cP(E_k)P(E_l)}{\sum_{i=2m_k+1}^{2m_l-1}\frac{1}{D_Y(i)}}
\leq \frac{cP(E_k)P(E_l)D_Y(2m_l)}{m_l-m_k}\no\\
&\le\frac{cP(E_k)P(E_l)F_Y(2m_l)}{m_l-m_k} \leq \frac{cP(E_k)}{l\log l-k\log k}.\no
\end{align}
As a consequence, we get
\begin{align}
\sum_{l=k+1}^{\ell_1-1}&P(E_kE_l)\leq cP(E_k)\sum_{l=k+1}^{\ell_1-1}\frac{1}{l\log l-k\log k}\nonumber\\
&\leq cP(E_k)\frac{1}{\log k}\sum_{l=k+1}^{\ell_1-1}\frac{1}{l-k}\leq cP(E_k)\frac{\log \ell_1}{\log k}\le cP(E_k),\label{e8}
\end{align}where for the last inequality, we use the fact
$
\frac{\log \ell_1}{\log k}\leq \gamma$
with some constant $\gamma$ depending only on $\varepsilon,$ whose proof is similar to \eqref{klg}.

Taking (\ref{e6}) and (\ref{e8}) together, we have
$$\sum_{k=k_3}^{N}\sum_{l=k+1}^{N}P(E_kE_l)\leq\sum_{k=k_3}^{N}\sum_{l=k+1}^{N}(1+\varepsilon)^2P(E_k)P(E_l)+\sum_{k=k_3}^{N}P(E_k).$$
Writing $H(\varepsilon)=(1+\varepsilon)^2$, owing to (\ref{syi}) , we have
\begin{align*}
\alpha_{H}&:=\liminf_{N\rightarrow\infty}\frac{\sum_{k=k_3}^{N}\sum_{l=k+1}^{N}P(E_kE_l)-\sum_{k=k_3}^{N}\sum_{l=k+1}^{N}HP(E_k)P(E_l)}{[\sum_{k=k_3}^{N}P(E_k)]^2}\nonumber\\
&\leq\lim_{N\rightarrow\infty}\frac{c}{\sum_{k=k_3}^{N}P(E_k)}=0.
\end{align*}
Applying  a generalized version of Borel-Cantelli lemma \cite[p. 235]{pe04}, we have
$$P(E_k,k\geq k_3 \text{ occur infinitely often})\geq\frac{1}{H+2\alpha_{H}}=\frac{1}{(1+\varepsilon)^2}.$$
Since $\varepsilon$ is arbitrary, letting $\varepsilon\rightarrow0,$ we conclude that
$$P(E_k,k\geq k_3 \text{ occur infinitely often})=1,$$ which implies \eqref{sfy}. The divergent part of Theorem \ref{main} for the chain $Y$ is proved. \qed

\subsubsection{Proof of Proposition \ref{p31}}

\proof Based on Proposition \ref{hp0}, the proof of Proposition \ref{p31} is a refinement of \cite[Proposition 3.1]{w19}. For convenience of the reader, we give the detail of its proof. To begin with, we prove (\ref{d1}). For $m\ge k\ge1,$ let $T_k,$ $h_k(j),\eta_{k,m}(j),j=1,2$ be as in \eqref{dtk}-\eqref{de2}.
If we write
\begin{align*}
&B_{k,1}:=\{Y_{T_k}=2k,2k+1 \text{ is a cutpoint}\},\\
&B_{k,2}:=\{Y_{T_k}=2k+1,2k \text{ is a cutpoint}\},
\end{align*}
then by the definition of cutpoints of the chain $Y,$ $\{k\in C^S\}=B_{k,1}\cup B_{k,2}$. Thus by the Markov property and \eqref{epcc}, it follows that
\begin{align*}
P(k\in C^S)&=P(B_{k,1})+P(B_{k,2})\\
&=h_k(1)\eta_{2k,2k}(2)Q_{2k+2}(2k+1,\infty,+)+h_k(2)Q_{2k+1}(2k,\infty,+)\\
&=h_k(1)\eta_{2k,2k}(2)\frac{1}{F_Y(2k+1)}+h_k(2)\frac{1}{F_Y(2k)},
\end{align*}
which leads to
\begin{align}\label{cul}
  \frac{h_k(1)\eta_{2k,2k}(2)+h_k(2)}{\max\{F_Y(2k),F_Y(2k+1)\}}\le P(k\in C^S)\le \frac{h_k(1)\eta_{2k,2k}(2)+h_k(2)}{\min\{F_Y(2k),F_Y(2k+1)\}}.
\end{align}
Noting that by assumption we have $\varrho_k^{-1}\rightarrow 1,$ as $k\rto,$ thus for any $\varepsilon>0,$  there exists a number $n_1$ such that for $n\ge n_1,$ $\varrho_n^{-1}>1-\varepsilon.$ Therefore for $n\ge n_1,$  we have $D_Y(n)=1+\sum_{j=n+1}^\infty \varrho_{n+1}^{-1}\cdots\varrho_j^{-1}\ge 1+\sum_{j=1}^{\infty}(1-\varepsilon)^j=\varepsilon^{-1}.$ Consequently, we have $D_Y(n)\rightarrow \infty$ as $n\rto.$ As a result, from \eqref{dyf}, we get $\lim_{n\rto}F_Y(n)=\infty.$
Since $\zeta_n\rightarrow \varrho^{-1}=1$ as $n\rto,$ then we have \begin{align}\label{dk}
  \lim_{k\rto}F_Y(k)/F_Y(k+1)=\lim_{k\rto}\z(1/F_Y(k+1)+\zeta_{k+1}\y)=1.
\end{align} Taking Proposition \ref{hp0} into account, from \eqref{cul} and \eqref{dk}, we obtain \eqref{d1}.

Next, we prove (\ref{d2}). Obviously, there is at most one cutpoint in $L_k.$
For $j<k,$ write \begin{align*}
&E_{jk}^{(1)}=\{\text{both }2j+1\text{ and }2k+1\text{ are cutpoints}\},\\
&E_{jk}^{(2)}=\{\text{both }2j+1\text{ and }2k\text{ are cutpoints}\},\\
&E_{jk}^{(3)}=\{\text{both }2j\text{ and }2k\text{ are cutpoints}\},\\
&E_{jk}^{(4)}=\{\text{both }2j\text{ and }2k+1\text{ are cutpoints}\}.
\end{align*}
Using the Markov property, we have
\begin{align*}
  P\left(E_{jk}^{(1)}\right)&=h_j(1)\eta_{2j,2j}(2)Q_{2j+2}^{2k}(2j+1,2k)\\
& \quad\quad\quad\quad \quad\quad  \times Q_{2k}^{2k+2}(2j+1,2k+1)Q_{2k+2}(2k+1,\infty,+).
\end{align*}
  The probabilities of $E_{jk}^{(i)}, i=2,3,4 $ can be computed similarly.
 We have
\begin{align}
P(j&\in C^S, k\in C^S)=\sum_{i=1}^4 P\left(E_{jk}^{(i)}\right)\nonumber\\
&=h_j(1)\eta_{2j,2j}(2)Q_{2j+2}^{2k}(2j+1,2k)\nonumber\\
&\quad\quad\quad\quad\quad\quad\times Q_{2k}^{2k+2}(2j+1,2k+1)Q_{2k+2}(2k+1,\infty,+)\nonumber\\
&\quad\quad+h_j(1){\eta_{2j,2j}(2)}Q_{2j+2}^{2k+1}(2j+1,2k)Q_{2k+1}(2k,\infty,+)\nonumber\\
&\quad\quad+h_j(2)Q_{2j+1}^{2k+1}(2j,2k)Q_{2k+1}(2k,\infty,+)\nonumber\\
&\quad\quad+h_j(2)Q_{2j+1}^{2k}(2j,2k){Q_{2k}^{2k+2}(2j,2k+1)}Q_{2k+2}(2k+1,\infty,+). \no
 \end{align}
Then, taking (\ref{epcb}), (\ref{epcc}), \eqref{cul} and  (\ref{dk}) into account, we get
 \begin{align}\label{jke}
   P(j&\in C^S, k\in C^S)\nonumber\\
   &\le \frac{1}{\min\{F_Y(2k),F_Y(2k+1)\}}\Big [h_j(1){\eta_{2j,2j}(2)}\nonumber\\
   &\quad\quad \times \z\{Q_{2j+2}^{2k}(2j+1,2k){Q_{2k}^{2k+2}(2j+1,2k+1)}+Q_{2j+2}^{2k+1}(2j+1,2k)\y\}\nonumber\\
&\quad\quad+ h_j(2)\z\{Q_{2j+1}^{2k}(2j,2k){Q_{2k}^{2k+2}(2j,2k+1)}+Q_{2j+1}^{2k+1}(2j,2k)\y\}\Big]\nonumber\\
 &\le \frac{ h_j(1){\eta_{2j,2j}(2)}+h_j(2)}{\min\{F_Y(2j),F_Y(2j+1)\}}\times \frac{h_k(1){\eta_{2k,2k}(2)}+h_k(2)}{\min\{F_Y(2k),F_Y(2k+1)\}}\nonumber\\
 & \quad\quad\times \max\left\{Q_{2j+2}(2j+1,2k,+),Q_{2j+1}(2j,2k,+)\right\}\times\frac{\max\left\{A(j,k),B(j,k)\right\}}{h_k(1)\eta_{2k,2k}(2)+h_k(2)}\nonumber\\
&\quad\quad\times\min\{F_Y(2j),F_Y(2j+1)\}\no\\
 &\le P(j\in C^S)P(k\in C^S)\max\left\{\frac{F_Y(2j+1)}{F_Y(2j+1,2k)},\frac{F_Y(2j)}{F_Y(2j,2k)}\right\}\nonumber\\
 &\quad\quad\quad\quad\quad\quad\quad\quad\quad\quad\quad\quad\times\frac{\max\left\{A(j,k),B(j,k)\right\}}{h_k(1)\eta_{2k,2k}(2)+h_k(2)},
  \end{align}
where
\begin{equation*}
\label{aij}A(j,k)=\frac{Q_{2j+2}^{2k}(2j+1,2k)}{Q_{2j+2}(2j+1,2k,+)}{Q_{2k}^{2k+2}(2j+1,2k+1)}+\frac{Q_{2j+2}^{2k+1}(2j+1,2k)}{Q_{2j+2}(2j+1,2k,+)}
\end{equation*}
 and \begin{equation*}\label{bij}B(j,k)=\frac{Q_{2j+1}^{2k}(2j,2k)}{Q_{2j+1}(2j,2k,+)}{Q_{2k}^{2k+2}(2j,2k+1)}+\frac{Q_{2j+1}^{2k+1}(2j,2k)}{Q_{2j+1}(2j,2k,+)}.\end{equation*}
 Now fix $\varepsilon>0.$ It follows from \eqref{lhe} and \eqref{ql} that \begin{align}
   \lim_{k\rto}A(j,k)=\lim_{k\rto}B(j,k)=\lim_{k\rto}h_k(1)\eta_{2k,2k}(2)+h_k(2)=-\frac{2\sigma}{1-\sigma},\no
 \end{align}
 and $\exists N_2>0 $ such that for all $k-j>N_2,$
 \begin{align}\label{ab}
   \frac{\max\left\{A(j,k),B(j,k)\right\}}{h_k(1)\eta_{2k,2k}(2)+h_k(2)}<\sqrt{1+\varepsilon}.
 \end{align}
On the other hand, by the definition of $F_Y(m,n)$ it is easy to see that \begin{align}
  \frac{F_Y(2j,2k)}{F_Y(2j+1,2k)}=\frac{1}{F_Y(2j+1, 2k)} +\zeta_{2j+1}.\label{fyf}
\end{align}
Let $\eta>0$ be an arbitrary small number.  Since $\varrho=1,$ in view of \eqref{ltz}, there exists $k_4>0$ such that for all $j\ge k_4,$ $1-\eta<\zeta_j<1+\eta.$
Thus, for $k> j>k_4,$ we have
\begin{align}
  F_Y(2j+1,2k)\ge 1+\sum_{i=1}^{2k-2j-2}(1-\eta)^i=\frac{1-(1-\eta)^{2k-2j-1}}{\eta}.\no
\end{align}
Therefore, we can find a number $N_3>0$ such that for all $k-N_3>j\ge k_4,$ $\frac{1}{F_Y(2j+1,2k)}<\eta.$
As a result, from \eqref{fyf} we see that
\begin{align}\label{fymf}
  \frac{F_Y(2j,2k)}{F_Y(2j+1,2k)}<1+2\eta, \text{ for all }k-N_3>j>k_4.
\end{align}
Taking \eqref{dk} and \eqref{fymf} together, we deduce that there exists a number $k_0$ such that for all $k-N_3>j>k_0,$
\begin{align}\label{mf}
  \max\left\{\frac{F_Y(2j+1)}{F_Y(2j+1,2k)},\frac{F_Y(2j)}{F_Y(2j,2k)}\right\}\le \sqrt{1+\varepsilon}\frac{F_Y(2j+1)}{F_Y(2j+1,2k)}.
\end{align}
Let $N=N_2\vee N_3.$  Substituting (\ref{ab}) and \eqref{mf}  into the rightmost hand of (\ref{jke}), we conclude that for $k-N>j>k_0,$
 \begin{align*}
   P(j\in C^S, k\in C^S)\le (1+\varepsilon)P(j\in C^S)P(k\in C^S)\frac{F_Y(2j+1)}{F_Y(2j+1,2k)}
 \end{align*}
  which completes the proof of (\ref{d2}).  \qed

\section{Number of cutpoints in $[2,n]$-Proof of Theorem \ref{nc}}\label{pn}

\proof We prove only the case $Z=Y,$ since the case $Z=X$ is easier. Suppose that $n$ is an odd number and $n=2m+1$ for some $m\ge1.$ Recall that $L_k=\{2k,2k+1\}, k\ge1.$ For $k\ge1,$ set $\eta_k=\z\{\begin{array}{cc}
              1, & k\in C^S \\
              0 &  k\notin C^S
            \end{array}
\y..$ Then we have $S_{n}=\sum_{k=1}^m \eta_k.$ It follows from Proposition \ref{p31} that $P(\eta_k=1)=P(k\in C^S)\sim \frac{c}{F_Y(2k)}, k\rto.$ Thus applying Lemma \ref{adf}, we see that for some proper constant $0<c_{10}<c_{11}<\infty,$  $\frac{c_{10}}{D_Y(2k)}<P(\eta_k=1)<\frac{c_{11}}{D_Y(2k)}, k\ge1.$ Consequently, we get $c_{10}\sum_{k=1}^m \frac{1}{D_Y(2k)}\le ES_n \le c_{11}\sum_{k=1}^m \frac{1}{D_Y(2k)},k\ge1.$ Therefore, on accounting of \eqref{dccy}, we have
\begin{align}\label{snb}
  c_{12}\sum_{k=2}^m \frac{1}{2k(\log\log(2k))^\beta}\le ES_n \le c_{13}\sum_{k=2}^m \frac{1}{2k(\log\log(2k))^\beta}
\end{align}
for some constants $0<c_{12}<c_{13}<\infty.$
 Notice that \begin{align*}
  \sum_{k=2}^m &\frac{1}{2k(\log\log(2k))^\beta}\sim \int_{2}^m\frac{dx }{2x(\log\log(2x))^{\beta}}\sim \frac{1}{2}\frac{1}{\log m(\log\log m)^{-\beta}}\no\\
  &\sim \frac{1}{2}\frac{1}{\log(2m+1)(\log\log (2m+1))^{-\beta}}= \frac{1}{2}\frac{1}{\log n(\log\log n)^{-\beta}},n\rto.
  \end{align*}
Thus from \eqref{snb} we get \eqref{lsn}.

In order to prove \eqref{asl}, noticing that
 $S_n$ is nonnegative and nondecreasing in $n$, thus on accounting of \eqref{snb} we have
\begin{align}
E\z(\max_{2\leq k\leq n}S_k\y)=E(S_n)\leq c_{13}\sum_{k=2}^m \frac{1}{2k(\log\log(2k))^\beta}\le c_{13}\sum_{k=3}^n\frac{1}{k(\log\log k)^\beta}.\label{snn1}
\end{align}
But for each $\varepsilon> 0$,
\begin{align}
\sum_{k=1}^{n}\frac{1}{k(\log{\log k})^{\beta}(\log{k})^{1+\varepsilon}(\log{\log k})^{-\beta}}<\infty.\label{snn2}
\end{align}
Therefore, with (\ref{snn1}) and (\ref{snn2}) in hands, applying \cite[Theorem 2.1]{ww20}, we get \eqref{asl}.
The theorem is proved for the chain $Y.$ \qed

\vspace{.5cm}

\noindent{{\bf \Large Acknowledgements:}} This project is partially supported by National
Natural Science Foundation of China (Grant No. 12071003; 11501008) and Nature Science Foundation
of Anhui Educational Committee (Grant No. YJS20210176).

\end{document}